\documentclass[final]{siamltex}         
\usepackage{graphicx,bm,amssymb,amsmath,xcolor} 
\newcommand{\bi}{\begin{itemize}} 
\newcommand{\ei}{\end{itemize}}
\newcommand{\ben}{\begin{enumerate}}
\newcommand{\een}{\end{enumerate}}
\newcommand{\be}{\begin{equation}}
\newcommand{\ee}{\end{equation}}
\newcommand{\bea}{\begin{eqnarray}} 
\newcommand{\eea}{\end{eqnarray}}
\newcommand{\ba}{\begin{align}} 
\newcommand{\ea}{\end{align}}
\newcommand{\bse}{\begin{subequations}} 
\newcommand{\ese}{\end{subequations}}
\newcommand{\bc}{\begin{center}}
\newcommand{\ec}{\end{center}}
\newcommand{\bfi}{\begin{figure}}
\newcommand{\efi}{\end{figure}}
\newcommand{\ca}[2]{\caption{#1 \label{#2}}}
\newcommand{\ig}[2]{\includegraphics[#1]{#2}}
\newcommand{\mbf}[1]{{\mathbf #1}}
\newcommand{\pO}{\Gamma} 
\newcommand{\Oc}{\Omega^c} 
\newcommand{\half}{\mbox{\small $\frac{1}{2}$}}

\newcommand{\RR}{\mathbb{R}^2}
\newcommand{\Dr}{{\cal D}} 
\newcommand{\Sr}{{\cal S}} 
\newcommand{\ttau}{\tilde\tau}                   
\DeclareMathOperator{\re}{Re}
\DeclareMathOperator{\im}{Im}
\newcommand{\matlab}{MATLAB}
\newtheorem{thm}{Theorem}

\newtheorem{rmk}[thm]{Remark}
\newcommand{\fref}[1]{Fig.~\ref{#1}}          
\newcommand{\sref}[1]{Sec.~\ref{#1}}          
\newcommand{\tref}[1]{Table~\ref{#1}}

\newcommand{\rref}[1]{Remark~\ref{#1}}
\newcommand{\lir}{\log\frac{1}{\rho}}     
\newcommand{\blir}{\biggl(\lir\biggr)}   
\newcommand{\Sop}{\mathbb{S}}         
\newcommand{\Dop}{\mathbb{D}}         

\newcommand{\Del}{\nabla}

\newcommand{\SLP}{\mathbf{S}}
\newcommand{\DLP}{\mathbf{D}}

\begin{document}

\title{Spectrally-accurate quadratures for evaluation of layer potentials
close to the boundary for the 
2D Stokes and Laplace equations
}

\author{Alex Barnett\thanks{Department of Mathematics, 
Dartmouth College, Hanover, NH, 03755, USA. {\em email:} {\tt ahb@math.dartmouth.edu}}%
\and Bowei Wu\thanks{Department of Mathematics, 
University of Michigan, Ann Arbor, MI, 48109, USA. {\em email:} {\tt boweiwu@umich.edu}.}
\and Shravan Veerapaneni\thanks{Department of Mathematics, 
University of Michigan, Ann Arbor, MI, 48109, USA. {\em email:} {\tt shravan@umich.edu}}}
\date{\today}
\maketitle

\begin{abstract}
%
Dense particulate flow simulations using integral
equation methods demand accurate evaluation of Stokes layer potentials
on arbitrarily close interfaces.
In this paper, we generalize
techniques for close evaluation of Laplace double-layer
potentials in
J.~Helsing and R.~Ojala, {\it J. Comput. Phys.} {\bf 227} (2008) 2899--2921.
We create a ``globally compensated'' trapezoid rule quadrature for the
Laplace single-layer potential on the interior and exterior of smooth curves.
This exploits a complex representation, 
a product quadrature (in the style of Kress) for the sawtooth function,
careful attention to branch cuts,
and second-kind 
barycentric-type formulae for Cauchy integrals and their derivatives.
%
Upon this we build accurate single- and double-layer Stokes potential evaluators
by expressing them in terms of Laplace potentials. 
We test their convergence
for vesicle-vesicle interactions, 
for an extensive set of
Laplace and Stokes problems,
and when
applying the system matrix in a boundary value problem solver
in the exterior of multiple close-to-touching ellipses.
We achieve typically 12 digits of accuracy
using very small numbers of discretization nodes per curve. 
We provide documented codes for other researchers to use.
\end{abstract}


\section{Introduction}
\begin{figure} 
\centering
\includegraphics[width=\textwidth]{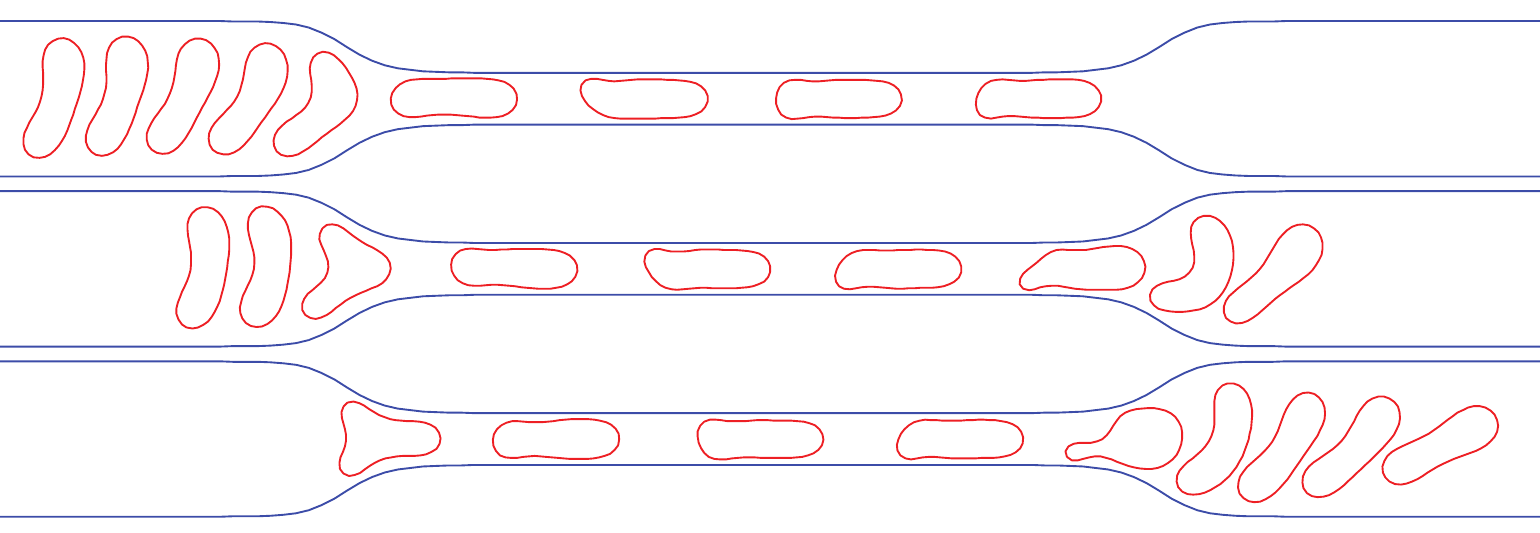} 
 \caption{Motivational example:
snapshots from a Stokes simulation of nine vesicles squeezing through a fixed-wall microfluidic device.
Using the close-evaluation scheme presented here,
accurate simulations can be run with only 32 points per vesicle.
(Without the scheme, instabilities creep in and the simulation breaks down after a few time-steps.)
Details of this application and its accuracy tests will be presented in
\cite{periodicp}.
\label{convDiv}
}
\end{figure}

Dense suspensions of deformable particles in viscous fluids are ubiquitous in natural and engineering systems. Examples include drop, bubble, vesicle, swimmer and blood cell suspensions. Unlike simple Newtonian fluids, the laws describing their flow behavior are not well established, owing to the complex interplay between the deformable microstructure and the macroscale flow. Besides experiments, direct numerical simulations are often the only means for gaining insights into the non-equilibrium behavior of such complex fluids. One of the main challenges for existing numerical methods to simulate {\em dense} or {\em concentrated} suspensions is to accurately resolve the particle-particle or the particle-wall interactions. In these complex fluid flows, more often than not, particles approach 
very
close to each other when subjected to flow (e.g., see Figure \ref{convDiv}). Numerical instabilities arise when the near interactions are not computed accurately, jeopardizing the entire simulation.

Boundary integral methods are particularly well-suited for vanishing Reynolds number problems where the Stokes equations govern the ambient fluid flow
\cite{pozrikidis}.
The advantages over grid- or mesh-based discretizations include:
a much smaller number of unknowns (exploiting the reduced
dimensionality), no need for smearing of interface forces onto a grid,
the availability of very high-order discretizations,
and of accelerated solvers such as the fast multipole method (FMM) \cite{lapFMM}
for handling the dense matrices in linear time.
%

The kernels of the integral operators such as the single-layer potential, 
\begin{equation} \label{S}
(\Sr \tau)(x) = \int_\pO G(x,y) \tau(y) \,ds_y
~,
\qquad x\in\RR, \end{equation}
defined on a smooth
closed planar curve $\Gamma$ for some smooth density function
$\tau$,
become {\em nearly singular} when the target point $x$ is close to $\Gamma$.
Neither smooth quadrature rules (such as the trapezoidal rule) nor singular quadratures are effective (uniformly convergent) for nearly singular integrals;
for example, the error in a fixed smooth quadrature rule grows
exponentially to ${\mathcal O}(1)$ as $x$ approaches $\Gamma$ \cite[Thm.~3]{ce}.
The objective of this paper is
the design of numerical integration schemes for \eqref{S} that,
given the smooth density $\tau$
sampled at the $N$ nodes of a periodic trapezoid rule on $\Gamma$,
exhibit {\em superalgebraic convergence} in $N$, with rate
{\em independent} of
the distance of $x$ from $\Gamma$.

\begin{rmk}
For smooth geometries and data many simple boundary integral solution methods,
such as Nystr\"om's method \cite{LIE},
exhibit superalgebraic convergence in the density.
Our goal is thus to provide layer potential evaluations
that are as accurate as the $N$-node spectral interpolant to the density itself,
i.e.\ {\em limited only by the data samples available}.
This means that in simulations, even those with close-to-touching
geometry, only the smallest number of unknowns $N$ required to
capture the density is needed, and optimal efficiency results.
\end{rmk}

Specifically, we develop a new suite of tools for evaluating layer potentials such as \eqref{S} (and their derivatives) on smooth closed curves for Laplace's equation in two dimensions (2D),
and from this build evaluators for Stokes potentials that can handle close-to-touching geometries and flow field evaluations arbitrarily close to curves,
with accuracies approaching machine precision.

Despite its importance in practical applications, very few studies
have addressed the accuracy issue with nearly singular integrals.
Adaptive quadrature on a target-by-target basis is impractically slow.
Beale et al.\ \cite{Beale,tlupova2013nearly} proposed a regularized
kernel approach that attains third-order accuracy by adding
analytically determined corrections.
%
Ying et al.\ \cite{Ying06} developed a method that
interpolates the potentials along extended surface
normals. Quaife--Biros \cite{quaife} applied this in the
context of 2D vesicle flows, attaining 5th-order accuracy.
Despite its high-order accuracy, the computational cost scales as a suboptimal
$\mathcal{O}(N^{3/2})$,
due to the need to compute values near the boundary using an upsampled
trapezoid rule.
Helsing--Ojala \cite{helsing_close}, exploiting Cauchy's theorem
and recurrence relations,
developed 16th-order panel-based close evaluation schemes in 2D;
these have recently been adapted to the complex biharmonic
formulation of Stokes potentials
by Ojala--Tornberg \cite{ojalastokes}.
The recent QBX scheme \cite{ce,qbx}
can achieve arbitrarily high order for Laplace and Helmholtz
potentials in 2D and 3D, but this
requires upsampling the density by a factor of 4--6.

The pioneering work of Helsing--Ojala includes a
``globally compensated'' scheme for the 2D double-layer potential (DLP)
\cite[Sec.~3]{helsing_close} which builds upon a second-kind
barycentric-type formula for quadrature of Cauchy's theorem
due to Ioakimidis et al.\ \cite{ioak}.
The scheme we present extends this to
the interior and exterior 2D single-layer potential (SLP);
since the complex logarithmic kernel is not single-valued,
this requires careful application of a spectrally-accurate
product quadrature for
the sawtooth function, in the style of Kress \cite{kress91}.
We also supply a true barycentric evaluation for
first derivatives of layer potentials that is stable for target
points arbitrarily close to nodes.
Unlike \cite{helsing_close}, we prefer to use an underlying global
quadrature (the periodic trapezoid rule) on $\Gamma$, since it is
most commonly used for vesicle simulations such as \fref{convDiv},
and (as our results show) is
somewhat more efficient in terms of $N$
than panel-based quadratures in this setting.%
\footnote{It is worth noting that even some panel-based schemes exploit
global schemes for adaptivity \cite{ojalastokes}.}

One advantage of our approach is that no auxiliary nodes
or upsampling is needed; another is that
the resulting discrete Cauchy sums are amenable for fast
summation via the FMM. While the new
tools for the planar Laplace equation are of interest in their own
right, our main motivation and interest is to enable accurate close
evaluation of Stokes potentials to target applications in interfacial
fluid mechanics. We accomplish this using the well-known fact
that Stokes potentials can be written in terms of Laplace potentials
and their derivatives \cite{fu2000fast, tornberg2008fast,
  ves2d}. Specifically in 2D, the Stokes SLP requires three Laplace SLP
evaluations and the DLP requires five Laplace DLP evaluations
(two of which are Cauchy-type).
We demonstrate in several numerical experiments that
our Stokes evaluations are very nearly as accurate as for Laplace.

The paper is organized as follows.  We define Laplace and Stokes
integral representations and set up notation in Section \ref{s:lps}.
Barycentric-type formulae for Cauchy integrals and their derivatives in
the interior and exterior of $\Gamma$ are discussed in Section
\ref{s:cau}. Based on these formulae, we derive spectrally accurate
global quadratures for close evaluation of Laplace potentials in
Section \ref{s:lap}, and demonstrate their performance
in evaluating all four types of boundary value problem (BVP) solutions.
We test performance of the quadratures applied to Stokes
potentials, and for the four Stokes BVPs, including one with close-to-touching
boundaries, in Section \ref{s:stokes}.
Finally, we summarize and discuss future work in Section \ref{s:conclusion}.

\section{Laplace and Stokes layer potentials}
\label{s:lps}

Let $\pO$ be a smooth closed Jordan curve in $\RR$, with outwards-directed
unit normal $n_y$ at the point $y\in\pO$.
Let $\Omega$ be the interior domain of $\pO$, and
$\Oc:=\mathbb{R}^2\backslash\overline{\Omega}$ be the exterior domain.
Let $\tau \in C(\pO)$ be a density function.
We review some standard definitions \cite[Ch.~6]{LIE}.
The Laplace SLP is defined by
\be
(\Sr\tau)(x) := \frac{1}{2\pi} \int_\pO
\blir
\tau(y) \,ds_y
\qquad x\in\RR
~.
\label{slp}
\ee
where the distance is $\rho := |r|$,
the displacement $r := x-y$, 
and $|x|:=\sqrt{x_1^2+x_2^2}$ is the Euclidean length of $x\in\RR$.
Finally, $ds_y$ is the arc length element on $\pO$.
The Laplace DLP is defined by
\be
(\Dr\tau)(x) := \frac{1}{2\pi} \int_\pO
\biggl(\frac{\partial}{\partial n_y}\lir\biggr)
\tau(y) \,ds_y
= \frac{1}{2\pi} \int_\pO
\biggl(\frac{r \cdot n_y}{\rho^2}\biggr)
\tau(y) \,ds_y
\qquad x\in\RR\backslash\pO
~.
\label{dlp}
\ee
Associating $\mathbb{C}$ with $\mathbb{R}^2$,
and noticing that the complex line element is $dy = i n_y ds_y$,
for purely real $\tau$ the DLP may also be written
as the real part of a Cauchy integral, as follows,
\be
(\Dr\tau)(x) = \re v(x)~,
\qquad \mbox{ where } \qquad
v(x) := \frac{1}{2\pi i}\int_\pO \frac{\tau(y)}{x-y} dy
\qquad x\in\mathbb{C}\backslash\pO
~.
\label{dlpcau}
\ee

Let $S$ be the restriction of $\Sr$ to evaluation on $\pO$,
in other words $S$ is the boundary integral operator with
kernel $k(x,y) = (1/2\pi) \log 1/\rho$.
$(S\tau)(x)$ exists as an improper integral.
Let $D$ be the restriction of $\Dr$ to evaluation on $\pO$,
in other words $D$ is the boundary integral operator with
kernel $k(x,y) = (1/2\pi) (r\cdot n_y / \rho^2)$,
taken in the principal value sense \cite[Sec.~6.3]{LIE}.
$D$ has a smooth kernel when $\Gamma$ is smooth.

We define the interior and exterior boundary limits of a function $u$
defined in $\RR\backslash\pO$,
at the point $x\in\pO$,
by $u^\pm(x) := \lim_{h\to 0^+} u(x \pm hn_x)$.
Likewise, $u^\pm_n(x) := \lim_{h\to 0^+} n_x \cdot \nabla u(x \pm hn_x)$.
We will need the following standard jump relations \cite[Sec.~6.3]{LIE}.
For any $C^2$-smooth curve $\pO$, and density function $\tau \in C(\pO)$,
\bea
(\Sr\tau)_n^\pm & = & (D^T \mp \half)\tau
\label{sjr}
\\
(\Dr\tau)^\pm & = & (D \pm \half)\tau
\label{djr}
~.
\eea

\subsection{Stokes potentials expressed via Laplace potentials}
\label{StokesViaLaplace}

Let $\sigma(y) = (\sigma_1(y), \sigma_2(y))$, for $y \in \Gamma$, be a smooth real-valued vector density function. The Stokes single and double layer potentials, denoted by $\SLP$ and $\DLP$, are defined by 
\begin{eqnarray}
(\SLP \mathbf{\sigma})(x) &:=& {1\over4\pi}\int_{\Gamma}\left( \lir\, I + \dfrac{ r\otimes r}{\rho^2} \right) \sigma(y) ds_{y}, 
\label{single_layer_kernel} \\ [0.5em]
(\DLP \sigma)(x) &:=& {1\over\pi}\int_{\Gamma}\left(\dfrac{ r\cdot  n_{y}}{\rho^2}\dfrac{ r\otimes r}{\rho^2} \right) \sigma(y) ds_{y}, 
\label{double_layer_kernel} 
\end{eqnarray}
where $r := x-y$ and $\rho := |r|$.
In \cite{fu2000fast, tornberg2008fast, ves2d}, fast algorithms to compute Stokes potentials were developed by expressing them in terms of Laplace potentials for which standard FMMs are applicable. We use the same strategy in this paper for close evaluation of Stokes potentials. Using the identity
\[ \dfrac{ r\otimes r}{\rho^2}\sigma = \dfrac{ r}{\rho^2}( r\cdot\sigma) = ( r\cdot\sigma)\Del_{x}\log\rho, \]
we can rewrite the Stokes SLP in terms of the Laplace SLP \eqref{slp}
as
\begin{equation}
\begin{aligned}
(\SLP \sigma)(x) \;=\;& {1\over4\pi}\int_{\Gamma}\blir \, \sigma ds_{y} + {1\over4\pi}\Del\int_{\Gamma}\blir \, (y\cdot\sigma) ds_{y}\\
 & - {1\over4\pi}x_{1}\Del\int_{\Gamma}\blir \, \sigma_{1}ds_{y}  - {1\over4\pi}x_{2}\Del\int_{\Gamma}\blir \, \sigma_{2}ds_{y}
~,
\end{aligned}
 \label{slp_stokes_via_lps}
\end{equation} 
where $\nabla = \nabla_x$ is assumed from now on.
Therefore, three Laplace potentials (and their first derivatives),
with density functions $y\cdot\sigma$, $\sigma_1$, and $\sigma_2$,
need to be computed to evaluate the Stokes SLP. Similarly, using the identity
$$
\Del\left ( \dfrac{ r\cdot n_{y}}{\rho^2} \right )
=\dfrac{ n_{y}}{\rho^2} - ( r\cdot n_{y})\dfrac{2 r}{\rho^4}
~,
$$
the Stokes DLP \eqref{double_layer_kernel} can be written as
\begin{equation}
\begin{aligned}
(\DLP \sigma)(x)\;
=\;& {1\over2\pi}\int_{\Gamma}{ n_{y}\over\rho^2}( r\cdot\sigma)ds_{y}  + {1\over2\pi}\Del\int_{\Gamma}{ r\cdot n_{y}\over\rho^2}( y\cdot\sigma)ds_{y}\\
 & - {1\over2\pi}x_{1}\Del\int_{\Gamma}{ r\cdot n_{y}\over\rho^2}\sigma_{1}ds_{y}  - {1\over2\pi}x_{2}\Del\int_{\Gamma}{r\cdot n_{y}\over\rho^2}\sigma_{2}ds_{y}
~.
\end{aligned}
\label{dlp_stokes_via_lps}
\end{equation}
The last three terms require Laplace DLP \eqref{dlp}
potentials (and first derivatives) for the same three densities
$y\cdot\sigma$, $\sigma_1$, and $\sigma_2$.
However, the first term is not a strict Laplace DLP of the form \eqref{dlp}:
the derivative is taken in the $\sigma$ rather than normal $n_y$ direction.
Yet, it can fit into our framework via two DLPs
if we generalize slightly the
Cauchy expression for the DLP \eqref{dlpcau} to allow complex densities $\tau$,
thus, using vector notation for the two components,
\be
{1\over2\pi}\int_{\Gamma}{ n_{y}\over\rho^2}( r\cdot\sigma)ds_y
= 
\re {1\over2\pi i} \int_\Gamma \frac{(\tau_1,\tau_2)}{x-y} dy~,
\quad\mbox{ where } \;
\tau_1 = (\sigma_1+i\sigma_2)\frac{\re n_y}{n_y}~,
\;
\tau_2 = (\sigma_1+i\sigma_2)\frac{\im n_y}{n_y}~.
\label{dlp_stokes_bowei_term}
\ee
So, in total five Laplace DLPs are needed.
Equations \eqref{slp_stokes_via_lps}--\eqref{dlp_stokes_bowei_term}
allow us to compute Stokes potentials by simply applying accurate
(Cauchy-form) Laplace close evaluation schemes, which are the focus of the
next two sections.

\section{Barycentric approximation of the interior and exterior
Cauchy integral formulae}
\label{s:cau}

In this section we describe an efficient and accurate
method to approximate a holomorphic
function $v$ from its boundary data sampled on a set of quadrature nodes
on the closed curve $\pO$.
The interior case is review
of Ioakimidis et al.\ \cite{ioak}, but we extend the method
to the exterior case
in a different manner from Helsing--Ojala \cite[Sec.~3]{helsing_close},
and correct an accuracy problem in the standard formula for the first
derivative. 
In Sec.~\ref{s:cautest} we present results showing uniform accuracy
close to machine precision.
We associate $\RR$ with $\mathbb{C}$.

We fix a quadrature scheme on $\pO$,
namely a set of nodes $y_j \in\pO$,
$j=1,\ldots,N$, and corresponding weights $w_j$, $j=1,\ldots,N$,
such that
$$
\int_\pO f(y) ds_y \;\approx\; \sum_{j=1}^N w_j f(y_j)
$$
holds to high accuracy for all smooth enough functions $f$.
Let $\pO$ be parametrized by the $2\pi$-periodic
map $Z:[0,2\pi) \to \RR$, with $Z(t) = Z_1(t) + iZ_2(t)$, such that
$\pO = Z([0,2\pi))$,
and with ``speed'' $|Z'(t)|>0$ for all $0\le t<2\pi$.
Then probably the simplest global quadrature arises from
the $N$-point periodic trapezoid rule \cite{davisrabin}
with equal weights $2\pi/N$ and nodes
\be
s_j:=\frac{2\pi j}{N}, \qquad j=1,\ldots,N
~.
\label{sj}
\ee
By changing variable to arc-length on $\pO$, we get a boundary quadrature
\be
y_j = Z(s_j), \quad w_j = \frac{2\pi}{N}|Z'(s_j)|, \qquad j=1,\ldots,N
~.
\label{ptr}
\ee
It is well known that this rule can be exceptionally accurate:
since the periodic trapezoid rule is exponentially convergent in $N$
for analytic $2\pi$-periodic integrands \cite{davis59} \cite[Thm.~12.6]{LIE}
\cite{PTRtref},
the rule \eqref{ptr} is exponentially convergent 
when $Z_1$ and $Z_2$ are analytic (hence $\pO$ is an analytic curve),
and the integrand $f$ is analytic.
The exponential rate is controlled by the size of the region of analyticity.
In the merely smooth case we have superalgebraic convergence.

\subsection{Interior case}
\label{s:cauint}

Cauchy's formula states that
any function $v$ holomorphic in $\Omega$ whose limit
approaching $\pO$ from the inside is $v^-\in C(\pO)$
may be reconstructed from its boundary data alone:
\be
\frac{1}{2\pi i}\int_\pO \frac{v^-(y)}{y-x} dy
\;=\;
\left\{\begin{array}{ll}v(x),&x\in\Omega\\0, & x\in \Oc\end{array}\right.
\label{cau}
\ee
Note that we have taken care to specify the data as the inside limit $v^-$;
this matters later when $v$ will be {\em itself}
generated by a Cauchy integral.

By combining \eqref{cau} with the special
case $(1/2\pi i)\int_\pO 1/(y-x) dy = 1 $ for $x\in\Omega$, we have
\be
\int_\pO \frac{v^-(y) - v(x)}{y-x} dy = 0, \qquad x\in\Omega
~.
\label{anal}
\ee
Even as the target point $x$ approaches $\pO$, the integrand remains smooth
(e.g.\ the neighborhood in which it is analytic remains large)
because of the cancellation of the pole,
and hence the quadrature rule \eqref{ptr} is accurate.
Thus
$$\sum_{j=1}^N \frac{v^-_j - v(x)}{y_j-x} w_j \;\approx\; 0~.$$
Rearranging derives a way to approximate $v(x)$, given the vector
of values $v^-_j:=v^-(y_j)$, namely
\be
v(x) \;\approx\;
\left\{\begin{array}{ll}\displaystyle \frac{\sum_{j=1}^N \frac{v^-_j}{y_j-x} w_j}{\sum_{j=1}^N \frac{1}{y_j-x} w_j},& x\in\overline{\Omega}, \;x\neq y_i, \;i=1,\dots,N\\
v^-_i,& x=y_i\end{array}\right.
\label{vcomp}
\ee
We have included $x$ in the closure of $\Omega$ because
in practical settings with roundoff error, targets may fall on $\pO$.
The second formula is needed when $x$ hits a node.
We believe \eqref{vcomp} is due to Ioakimidis et al.\ \cite[(2.8)]{ioak};
Helsing--Ojala call it ``globally compensated'' quadrature
\cite{helsing_close}.
It is in fact a {\em barycentric} Lagrange polynomial interpolation formula
of the second form,
with two crucial differences from the usual setting \cite{berruttref}:
\bi
\item the nodes are no longer on the real axis, and
\item the weights $w_j$ come simply from quadrature weights on the curve
rather than from the usual formula related to Lagrange polynomials.
\ei
In the case of $\Omega$ the unit disc with equispaced nodes,
the equivalence of \eqref{vcomp} to barycentric
interpolation was recently explained by
Austin--Kravanja--Trefethen \cite[Sec.~2.6]{austinrou}.

A celebrated key feature of barycentric formulae is
numerical stability even as the evaluation point
$x$ approaches arbitrarily close to a node $y_i$.
Although relative error grows without limit in both numerator and denominator
of \eqref{vcomp}, due to roundoff error in the dominant terms $1/(y_i-x)$,
these errors cancel (see \cite[Sec.~7]{berruttref}, \cite{baryhigham}
and references within).
Thus close to full machine precision is attainable even in this limit.

\subsection{Exterior case}
\label{s:cauext}

We turn to the exterior Cauchy formula that states for $v$
holomorphic in $\Oc$,
\be
\frac{1}{2\pi i}\int_\pO \frac{v^+(y)}{y-x} dy \;=\; 
\left\{\begin{array}{ll}v_\infty,&x\in\Omega\\
v_\infty - v(x),& x\in \Oc
\end{array}\right.
\label{ve}
\ee
where $v_\infty := \lim_{x\to\infty} v(x)$.
We are interested only in $v$ that can be generated by the
exterior Cauchy integral, i.e.\ the case $v_\infty=0$.
We pick a simple non-vanishing function $p$ with $p_\infty=0$
that will play the role that the constant function played in the interior case.
We choose $p(x) = 1/(x-a)$ where $a\in\Omega$ is a fixed arbitrary
point chosen not near $\pO$.
Applying \eqref{ve} this is generated by
$$
\frac{1}{x-a} =  \frac{-1}{2\pi i}\int_\pO \frac{(y-a)^{-1}}{y-x} dy
$$
Multiplying both sides by $(x-a)$ we get a way to represent
the constant function 1 via an exterior representation.
Using this we create an exterior equivalent of \eqref{anal},
\be
\int_\pO \frac{v^+(y) - (y-a)^{-1}(x-a)v(x)}{y-x} dy = 0
\qquad \mbox{ for } x\in\Oc
~.
\label{anale}
\ee
The integrand remains smooth and analytic even as the target
$x$ approaches $\pO$.
Substituting the periodic trapezoid rule in \eqref{anale}
and rearranging as in the interior case gives
\be
v(x) \;\approx\;
\left\{\begin{array}{ll}\displaystyle \frac{1}{x-a} \cdot \frac{\sum_{j=1}^N \frac{v^+_j}{y_j-x} w_j}{\sum_{j=1}^N \frac{(y_j-a)^{-1}}{y_j-x} w_j},
&x \in \Oc,\;x\neq y_i, \;i=1,\dots,N\\
v^+_i,& x=y_i\end{array}\right.
\label{vcompe}
\ee
which is our formula for accurate evaluation of the
exterior Cauchy integral. It also has barycentric stability near nodes.

\begin{rmk} 
Helsing--Ojala \cite[Eq.~(27)]{helsing_close}
mention a different formula for the exterior case,
$$
v(x) \;\approx\;
\frac{\sum_{j=1}^N \frac{v^+_j}{y_j-x} w_j}{-2\pi i + \sum_{j=1}^N \frac{1}{y_j-x} w_j}
~, \qquad x \in \Oc,\;x\neq y_i, \;i=1,\dots,N.
$$
This is marginally simpler than our \eqref{vcompe} since the interior
point $a$ is not needed;
we have not compared the two methods numerically,
since our scheme performs so well.
\label{r:2pii}
\end{rmk}  

\subsection{First derivative and its barycentric form, interior case}
\label{s:caupint}

Accurate Stokes evaluation demands accurate first derivatives of Laplace
potentials and hence of the Cauchy representation.
The interior Cauchy formula for the first derivative is,
$$
v'(x) = \frac{1}{2\pi i}\int_\pO \frac{v^-(y)}{(y-x)^2} dy, 
\qquad x \in\Omega~.
$$
We can combine this with the Cauchy formula as in \eqref{anal} to get,
\be
\int_\pO \frac{v^-(y) - v(x) - (y-x) v'(x)}{(y-x)^2} dy = 0
\qquad x \in\Omega~,
\label{analp}
\ee
which holds because the middle term vanishes (the contour integral
of $1/(x-y)^n$ is zero for integer $n\neq 1$).
The integrand is analytic and smooth even as $x$ approaches $\pO$
because the numerator kills
Taylor terms zero and one in the expansion of $v$ about $x$,
so the trapezoid rule \eqref{ptr}, as before, is accurate.
Applying the quadrature (making sure to keep the middle term, which is
mathematically zero but numerically necessary to compensate the quadrature)
gives
\be
v'(x) \;\approx\;
\left\{\begin{array}{ll}\displaystyle
\frac{\sum_{j=1}^N \frac{v^-_j-v(x)}{(y_j-x)^2} w_j}{\sum_{j=1}^N \frac{1}{y_j-x} w_j},& x\in\overline{\Omega}, \;x\neq y_i, \;i=1,\dots,N\\
\displaystyle -\frac{1}{w_i}\sum_{j\neq i}\frac{v^-_j-v^-_i}{y_j-y_i}w_j,& x=y_i\end{array}\right.
\label{vpcomp}
\ee
Here the case where $x$ coincides with a node
is derived by taking the limit as $x$ approaches a node.
\eqref{vpcomp} is analogous to the derivative of the barycentric interpolant
derived by Schneider--Werner \cite[Prop.~11]{schneider86}
(who also generalized to higher derivatives).
The case for $x$ not a node is equivalent
to the derivative formula of Ioakimidis et al.\ \cite[(2.9)]{ioak},
and that of \cite[Algorithm {\tt P'}]{austinrou}.

However, \eqref{vpcomp} (and its
equivalent forms cited above) suffers from catastrophic cancellation
as $x$ approaches a node, a point
that we have not seen discussed in the literature.
Even if $v(x)$ is computed to high accuracy (say, using \eqref{vcomp}),
the {\em relative} accuracy of the term $(v_j^--v(x))$ in \eqref{vpcomp}
deteriorates, in a way that
involves no explicit cancellation of poles as in the second barycentric form.
A solution which regains true barycentric stability is
to evaluate this term via
\be
v_j^--v(x) \;\approx\;
\frac{\sum_{k\neq j} \frac{v^-_j-v^-_k}{y_k-x} w_k}{\sum_{k=1}^N \frac{1}{y_k-x} w_k}
~.
\label{dvbary}
\ee
Because the numerator term $k=j$ is absent,
then as $x$ tends to node $y_j$ an overall factor of $(y_j-x)$ dominates
in a way that cancels (to high accuracy)
one power in the term $(y_j-x)^2$ in \eqref{vpcomp}.
Thus using \eqref{dvbary} for each term $v_j^--v(x)$ is
a true second barycentric form, which we believe is new.
The problem is that it increases the effort from ${\mathcal O}(N)$ to
${\mathcal O}(N^2)$ per target point $x$.
Our remedy is to realize that \eqref{dvbary} is only helpful for
small distances $|y_j-x|$: thus we use the value $v(x)$ from
\eqref{vcomp} unless $|y_j-x|<\delta$, in which case \eqref{dvbary} is used.
We choose $\delta=10^{-2}$, since in most settings
with $\Omega$ of size ${\mathcal O}(1)$
only a small fraction of targets lie closer than this to a node,
and at most around 2 digits are lost due to the loss of barycentric stability
for larger target-node distances.
Unless many nodes are spaced much closer than $\delta$, the method
remains ${\mathcal O}(N)$ per target point.

\begin{rmk}  
It has been pointed out%
\footnote{L. N. Trefethen, personal communication.} to us
that another way to alleviate the issue of the first formula in \eqref{vpcomp}
not itself being barycentric is to instead compute $v'(y_i)$ at each
of the nodes using the second formula in \eqref{vpcomp}, then to
use barycentric interpolation \eqref{vcomp} from these values.
I.e.\ one {\em interpolates the derivative instead of differentiating the
interpolant}.
We postpone comparing these two methods to future work, since
we note that we already achieve close to machine precision errors, uniformly.
\label{r:tref}
\end{rmk}  

\subsection{First derivative with a barycentric form, exterior case}
\label{s:caupext}

Combining the ideas of \eqref{anale} and \eqref{analp} we have
the identity
$$
\int_\pO \frac{v^+(y) - v(x) - (x-a)(y-a)^{-1}(y-x) v'(x)}{(y-x)^2} dy
=0~,
\qquad \mbox{ for } x\in\Oc
~,
$$
which, as with \eqref{analp}, has smooth analytic integrand
even as $x$ approaches $\pO$
because the first two Taylor terms are cancelled.
Inserting the quadrature rule gives the approximation
\be
v'(x) \approx \frac{1}{x-a}\cdot  \frac{\sum_{j=1}^N \frac{v^+_j-v(x)}{(y_j-x)^2} y'_j w_j}{\sum_{j=1}^N \frac{(y_j-a)^{-1}}{y_j-x} y'_j w_j}~,
\qquad x\in \RR \backslash \Omega, \;x\neq y_i, \;i=1,\dots,N\
\label{vpcompe}
\ee
which, as in the interior case, does not have barycentric stability
as $x$ approaches a node.
Unfortunately the formula analogous to \eqref{dvbary} which uses \eqref{vcompe}
to write $v(x)$ also fails to give stability, because the $k=j$ term no
longer vanishes and roundoff in this term dominates.
However, it is easy to check that the following form is mathematically
equivalent,
\be
v_j^+-v(x) \;\approx\;
\frac{1}{x-a}\biggl[
\frac{\sum_{k\neq j} \frac{v^+_j(y_j-a)(y_k-a)^{-1}-v^+_k}{y_k-x} w_k}
{\sum_{k=1}^N \frac{(y_k-a)^{-1}}{y_k-x} w_k}
- (y_j-x)v^+_j
\biggr]
~,
\label{dvbarye}
\ee
and does
give barycentric stability when inserted into \eqref{vpcompe},
because, as $x$ approaches $y_j$, the factor $y_j-x$ in both terms is explicit.
As with the interior case, we only use this when $|y_j-x|<\delta$.
This completes our recipes for interior/exterior Cauchy values and first
derivatives.

\bfi[t] 
\ig{height=1.3in}{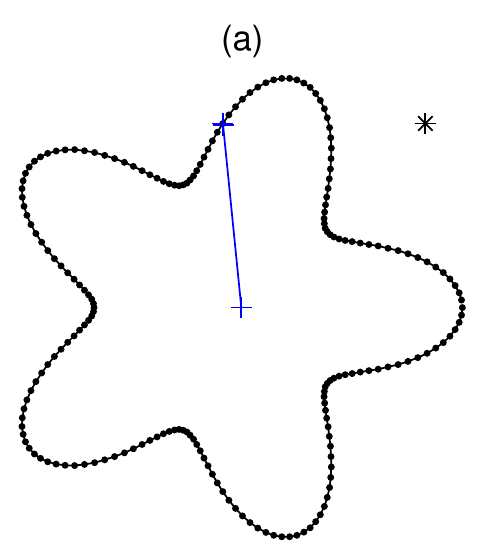} 
\ig{height=1.7in}{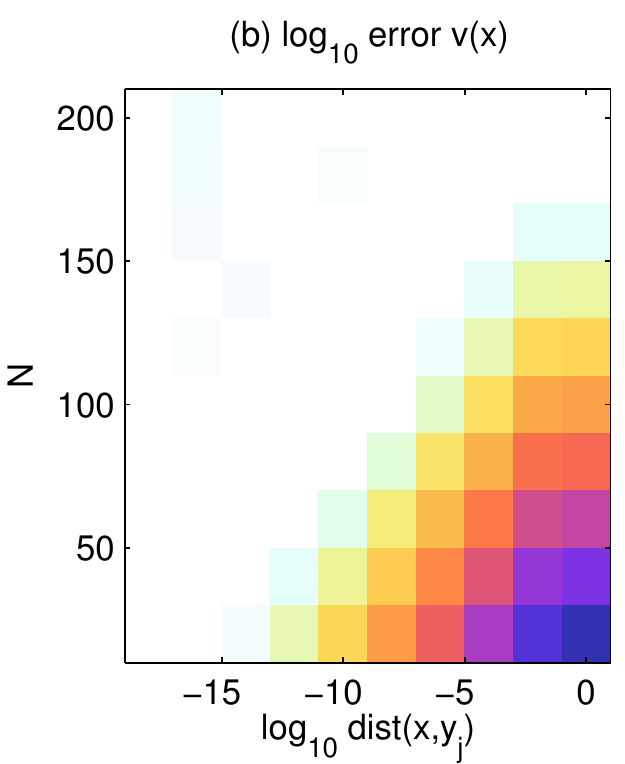}
\ig{height=1.7in}{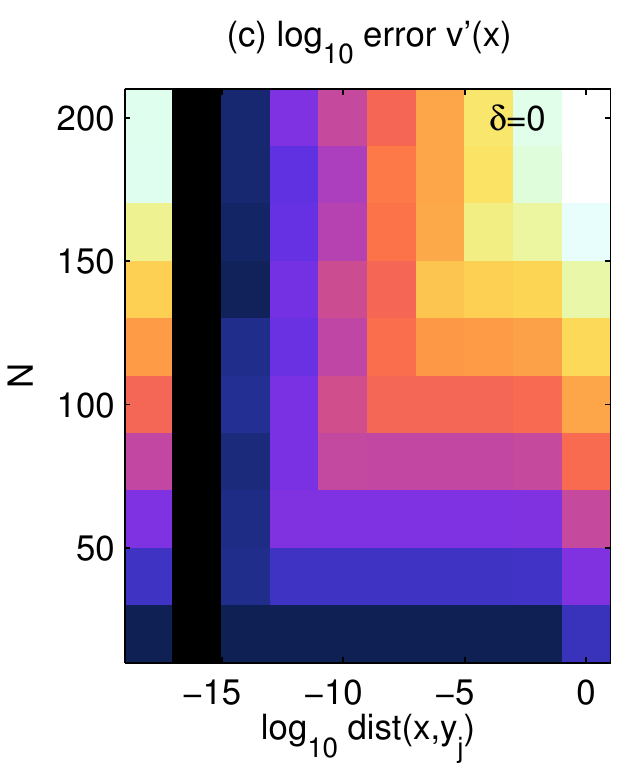}
\ig{height=1.7in}{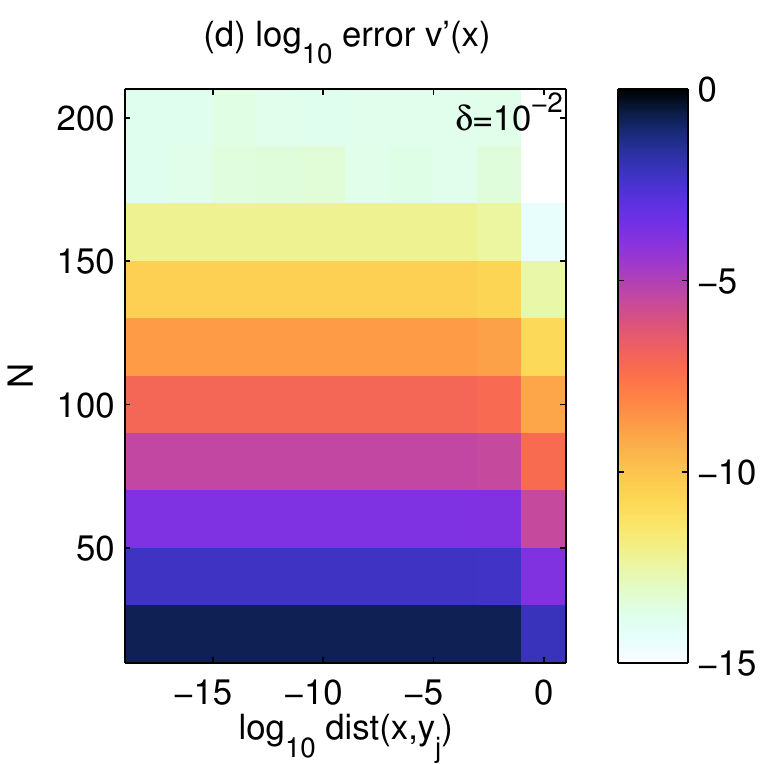}
\ca{
(a) Curve (with $N=200$ nodes shown)
for the test of the interior Cauchy integral evaluation;
points $x$ (shown by $+$ symbols along the straight line)
lie at approximate distances from a node
$0,10^{-16},10^{-14},\ldots,10^{-2},1$.
The test function is $v(x) = 1/(x-b)$ where the pole $b = 1.1+1i$
(shown by $\ast$) is a distance of 0.5 from $\pO$.
(b) Convergence of
$\log_{10}$ of the absolute error (see colorbar at right) for $v(x)$, using
formula \eqref{vcomp}, plotted vs distance of $x$ from
a node (horizontal axis) and number of nodes $N$ (vertical axis).
(c) Convergence of $v'(x)$, using the formula \eqref{vpcomp}.
(d) Convergence of $v'(x)$, inserting \eqref{dvbary} into \eqref{vpcomp}
for distances below $\delta=10^{-2}$.
Exterior results are very similar, so we do not show them.%
}{f:cauint}
\efi

\subsection{Numerical tests of values and derivatives close to the curve}
\label{s:cautest}

We test the Cauchy integral evaluation formulae presented in this section
on the smooth star-shaped domain of Fig.~\ref{f:cauint}(a)
given by the radial function $f(\theta) = 1 + 0.3 \cos(5\theta)$.
(This shape is also used in \cite{helsing_close}.)
For the interior case, the test points $x$ lie along the line shown
in Fig.~\ref{f:cauint}(a), at a set of distances from a node $y_j$
logarithmically spanning the range from machine precision to 1.
A test point $x=y_j$ is also included (data for this test point appears
in the left-most column of each plot (b)--(d)).
The function $v$ is a pole located outside of $\pO$ and is therefor
holomorphic in $\Omega$; its maximum magnitude on $\pO$ is of order 1.

Fig.~\ref{f:cauint}(b)
shows exponential convergence with 15 digit accuracy in value
reached by $N=180$
at all distances, and that small distances converge at the same rate
but with a smaller prefactor.
For comparison, applying the quadrature
scheme \eqref{ptr} directly to the Cauchy integral,
as is common practice, with the same $N=180$, gives 15 digits of accuracy
at the most distant point,
but zero digits of accuracy at all other points.

The derivative formula \eqref{vpcomp}
of 
Ioakimidis et al. is tested in plot (c):
there is a clear loss of accuracy in inverse proportion to the distance
from a node, regardless of $N$, simply
due to non-barycentric loss of relative error in the term $v^-_j-v(x)$.
(However, when $x$ coincides with a node, convergence
is again exponential).
Finally, evaluating $v^-_j-v(x)$ by \eqref{dvbary} for distances
less than $\delta=10^{-2}$ gives plot (d), which achieves 14 digit
accuracy at all distances.

We also tested the exterior methods,
choosing a generic interior point $a=-0.1$ in the method,
and using $v$ given by a pole at $b = 0.1+0.5i$ a distance 0.33 from $\pO$.
This $v$ is holomorphic in $\Oc$ and has $v_\infty =0$.
Results for the value formula \eqref{vcompe}, and
derivative formula \eqref{vpcompe} with $v^+_j-v(x)$ evaluated
via \eqref{dvbarye}, are essentially identical to the interior
case, with equally good achievable accuracies, so we do not show them.

Note that the {\em rate} of exponential convergence is
clearly affected by the choice of holomorphic test function $v$:
moving the pole of $v$ closer to $\pO$ worsens the rate since the data becomes
less smooth. We chose poles not too far from $\pO$.
For comparision, for an entire function, such as $v(x) = e^{2x}$,
full convergence in the interior is achieved at only $N=80$.

\section{Evaluation of Laplace layer potentials with global quadrature}
\label{s:lap}

The goal of this section is to describe accurate methods to evaluate
the single- and double-layer potential in the interior and exterior of
a closed curve $\pO$, given only the density values $\tau_j:=\tau(y_j)$
at the nodes $y_j$ belonging to a global quadrature \eqref{ptr} on the curve.
We remind the reader that
direct application of the rule \eqref{ptr} to the
layer potentials \eqref{slp} and \eqref{dlp} is highly inaccurate near $\pO$.
In contrast,
our methods retain accuracy and efficiency for targets $x$ arbitrarily
close to $\pO$. The effort will be ${\mathcal O}(N(N+M))$ for $N$ nodes and $M$ targets (although see discussion in \sref{s:conclusion}).

\subsection{Laplace double-layer potential}
\label{s:lapdlp}

Recall that the
double-layer potential \eqref{dlp} can be written as the real part of
the function $v$ given by the Cauchy integral \eqref{dlpcau}.
$v$ is holomorphic in $\Omega$ and in $\Oc$.
For interior evaluation, Helsing--Ojala \cite{helsing_close}
proposed a two-stage scheme, which for convenience we review
in our setting of the global periodic trapezoid quadrature:
\bi
\item [{\bf Step 1.}] Approximate the boundary data $v^-_j:=v^-(y_j)$ which
is the interior limit of the function \eqref{dlpcau} at each of the nodes.
\item [{\bf Step 2.}] Use this data to numerically approximate
the Cauchy integral \eqref{cau} to generate $v(x)$
at any $x\in\Omega$, using the method of Section~\ref{s:cauint}.
\ei
Finally, taking $u = \re v$ extracts the desired potential.
We extend this to first derivatives by including in Step 2 the
method of Section~\ref{s:caupint}
to evaluate $v'(x)$, then extracting the gradient as
$\nabla u = (\re v',-\im v')$.

It only remains to present Step 1 (following \cite[Sec.~3]{helsing_close}).
For $x\in\pO$, let $v^-(x) := \lim_{y\to x, y\in\Omega} v(y)$ be the
interior limit of \eqref{dlpcau}. It follows from the
Sokhotski--Plemelj jump relation \cite[Thm.~7.6]{LIE} that
$$
v^-(x) \;=\; -\half \tau(x) - \frac{1}{2\pi i} \mbox{PV}
\int_\pO \frac{\tau(y)}{y-x} dy~, \qquad x\in\pO
~,
$$
where PV indicates the principal value integral.
We split the PV integral into an analytic and Cauchy part:
$$\mbox{PV} \int_\pO \frac{\tau(y)}{y-x} dy \;=\;
\int_\pO \frac{\tau(y)-\tau(x)}{y-x} dy \;+\; \tau(x) \,\mbox{PV}
\int_\pO \frac{1}{y-x} dy~.$$
The latter integral is $-\half$ for $x\in\pO$ for any closed curve $\pO$.
Thus
\be
v^-(x) \;=\; -\tau(x) - \frac{1}{2\pi i}
\int_\pO \frac{\tau(y)-\tau(x)}{y-x} dy~, \qquad x\in\pO
\label{vm}
\ee
The integrand is analytic (compare \eqref{anal}),
so the periodic trapezoid rule is again excellent.
We need the boundary point $x=y_k$, i.e.\ the $k$th node, so must
use the correct limit
of the integrand at the diagonal $y=x$.
If we define $\ttau(t) := \tau(Z(t))$ as the density in the parameter
variable, and let $\tau'_k := \ttau'(s_k)$,
then the integrand has diagonal limit
$\ttau'(s_k)/|Z'(s_k)|$, and applying quadrature \eqref{ptr}
to \eqref{vm} gives
\be
v^-_k \;=\; -\tau_k - \frac{1}{2\pi i} \sum_{j \neq k} \frac{\tau_j-\tau_k}{y_j-y_k}
w_j - \frac{\tau'_k}{iN}~,
\qquad k=1,\ldots,N~.
\label{vmk}
\ee
We compute the vector $\{\tau'_j\}_{j=1}^N$ by spectral differentiation
\cite{tref}
of the vector $\{\tau_j\}_{j=1}^N$ via the $N$-point fast Fourier transform
(FFT). This completes Step 1 for the interior case.

The case of $x$ exterior to $\pO$ is almost identical, except that
by the Sokhotski--Plemelj jump relation we change Step 1 to
\be
v^+_k \;=\; v^-_k + \tau_k
\;=\;
- \frac{1}{2\pi i} \sum_{j \neq k} \frac{\tau_j-\tau_k}{y_j-y_k}
w_j - \frac{\tau'_k}{iN}~,
\qquad k=1,\ldots,N~,
\label{vpk}
\ee
and in Step 2 we now use
the exterior methods of Sections~\ref{s:cauext} and \ref{s:caupext}
for
$v(x)$ and $v'(x)$.

\subsection{Laplace single-layer potential}
\label{s:lapslp}

This section is the heart of the contribution of this paper.
Although we have not seen this in the literature (other than \cite{ce}),
it is also possible to
write the single-layer potential \eqref{slp} as the real part
of a holomorphic function,
\be
u(x) = (\Sr\tau)(x)
= \re v(x)
\ee
where the function $v$ is defined by
\be
v(x) := \frac{1}{2\pi}\int_\pO \left(\log\frac{1}{y-x}\right) \tau(y) \,|dy|
~,
\qquad x\in\RR\backslash\pO
~.
\label{slpc}
\ee
Note that the integration element $|dy| = dy/in_y$ is real.
We present the interior and exterior cases in turn.

\subsubsection{Interior case}

We first explain the (simpler) interior case $x\in\Omega$.
In order that $v$ be holomorphic in $\Omega$,
care needs to be taken with the branch cuts of the logarithm in \eqref{slpc},
and it must be considered as a function of two variables,
$L(y,x):=\log 1/(y-x)$, which agrees with the standard logarithm
up to the choice of Riemann sheet.
Fixing a boundary point $y_0\in\pO$,
the following branch choices for $L$ are sufficient for $v$ to be holomorphic
\cite[Remark~7]{ce}:
i) for each fixed $x\in\Omega$, as $y$ loops around $\pO$
a jump of $2\pi i$ occurs in $L$ only at $y_0$,
and ii) for each fixed $y\in\pO$, $L(y,\cdot)$ is continuous in $\Omega$.

As before we have two steps, only the first of which 
differs from the double-layer case:
Step 1 finds the interior boundary data $v^-_j$ for
the function \eqref{slpc},
then Step 2 evaluates $v$ at arbitrary interior target
points using \eqref{cau} and
the Cauchy integral method of Section~\ref{s:cauint}.
Finally the real part is taken.
Accurate first partials are also found in the same way as the double-layer case.

All that remains is to explain Step 1.
Inserting the parametrization of $\pO$ from Section~\ref{s:cau}
and recalling $\ttau(s) := \tau(Z(s))$,
the interior boundary data of $v$ at $x = Z(t)\in\pO$ is
\be
v^-(x) := \lim_{\Omega\ni z\to x} \frac{1}{2\pi} \int_{0}^{2\pi}
\left(\log\frac{1}{Z(s)-z}\right) \ttau(s) |Z'(s)| ds
~.
\label{vmslp}
\ee
To handle the logarithmically singular kernel which results in the limit,
we exploit the identity
\be
\log\frac{1}{Z(s)-Z(t)} \; = \; \log\frac{e^{is}-e^{it}}{Z(s)-Z(t)} - \log
(e^{is}-e^{it})
~.
\label{split}
\ee
For $\pO$ smooth,
the first term on the right-hand side is smooth as a function of $s$
when the correct branch cuts of $\log$ are taken,
while the second term
(relating to the single-layer kernel on the unit disc) can be handled
analytically.
Similar ideas are used in analytic PDE theory~\cite[Sec.~8]{tothzeld09}.
Inserting this into \eqref{vmslp}, realizing that the
sense of the limit (interior) only affects
the second term, gives
\be
v^-(Z(t)) := \frac{1}{2\pi} \int_{0}^{2\pi}
\left(\log\frac{e^{is}-e^{it}}{Z(s)-Z(t)}\right) \ttau(s) |Z'(s)| ds
-
\frac{1}{2\pi} \int_{0}^{2\pi}
\log(e^{is}-e^{i(t+i0)}) \ttau(s) |Z'(s)| ds
\label{vmspl}
\ee
where $t+i0$ indicates the limit where the imaginary part of
$t$ approaches zero from above, i.e.\ $e^{it}$ approaches the unit circle
from inside.
We write $\log(e^{is}-e^{i(t+i0)}) = is + g(s-t)$
where $g$ is a convolution kernel.
We simply drop the non-convolutional term $is$ because it contributes
a purely imaginary constant to $v$ that will have no effect on $\re v$.
The sense of the limit allows a Taylor expansion of the logarithm in the kernel,
\be
g(s) = \log(1 - e^{-i(s-i0)}) = \sum_{n>0} \frac{e^{-ins}}{n}
 =: \sum_{n\in\mathbb{Z}} g_n e^{ins}~,
\label{g}
\ee
where the last term defines a Fourier series
whose coefficients 
are read off as
\be
g_n = \left\{ \begin{array}{ll} -\frac{1}{n}, & n<0 \\ 0 ,& n\ge 0\end{array}
\right.
\label{gn}
\ee
These coefficients can be used to construct weights
that approximate the product quadrature with $g$,
via trigonometric polynomials.
Let any function $f$ be smooth and $2\pi$-periodic, then
$$
\int_0^{2\pi} f(s) g(s) ds \;\approx \; \sum_{j=1}^N R_j f(s_j)
$$
holds to spectral accuracy, where the general formula%
\footnote{We drop the factor of $\half$ from the last term
present in \cite{kress91} and \cite[Sec.~6]{hao}, it being
exponentially small.}
for the weights
(as derived in the review \cite[Sec.~6]{hao}),
followed by the weights for our particular case \eqref{g} of $g$, is
\be
R_j = \frac{2\pi}{N}\sum_{|n|<N/2} \overline{g}_n e^{-ins_j}
= -\frac{2\pi}{N}\sum_{n=-N}^{-1} \frac{e^{-ins_j}}{n}
~,
\qquad j=1,\ldots, N~.
\label{Rj}
\ee
The action of $g$ as a convolution kernel on any $f$ is thus approximated
by a circulant matrix:
$$
\int_0^{2\pi} f(s) g(s-s_k) ds \;\approx \; \sum_{j=1}^N R_{j-k} f(s_j),
\qquad k=1,\ldots, N
~.
$$

We are now ready to apply quadrature rules to \eqref{vmspl}.
For each boundary target $x=y_k=Z(s_k)$, the first (smooth) term
is well approximated using the
periodic trapezoid rule with the correct diagonal limit,
while the second term uses the above product quadrature, giving
\be
v^-_k \approx
\frac{1}{2\pi} \sum_{j\neq k} \left( \log \frac{e^{i(s_k-s_j)}}{Z(s_k)-Z(s_j)}
\right) w_j \tau_j + \frac{1}{2\pi} \log \frac{i e^{is_k}}{Z'(s_k)} w_k \tau_k
- \frac{1}{2\pi} \sum_{j=1}^N R_{j-k} w_j \tau_j
,\quad k=1,\ldots N.
\label{vmnyst}
\ee
With \eqref{vmnyst} and \eqref{Rj} defined, Step 1 is complete.

A couple of practical words are needed.
The vector $\{R_j\}_{j=1}^N$ is simply filled by taking the FFT of \eqref{gn}.
Handling the branch cuts of the first term in \eqref{vmnyst}
so that it corresponds to a smooth kernel
may seem daunting. In fact this is done easily by
filling the $N$-by-$N$ matrix {\tt S} corresponding to the first two terms in
\eqref{vmnyst}, using the machine's standard branch cut for log,
then applying to it the following simple \matlab\ code,
\vspace{.5ex}
\begin{verbatim}
for i=1:numel(S)-1
  p = imag(S(i+1)-S(i));
  S(i+1) = S(i+1) - 2i*pi*round(p/(2*pi));
end
\end{verbatim}
\vspace{.5ex}
This loops through all matrix elements
and adjusts them by an integer multiple of $2\pi i$ whenever they
jump by more than $\pi$ in imaginary part from a neighboring element.
Once $N$ is sufficiently large to resolve the kernel, such large jumps
cannot occur unless a branch is being crossed.

\subsubsection{Exterior case}
\label{s:lapslpext}

Here we describe the differences from the interior case.
In Step 1, the crucial sign change causes the convolutional part of the kernel
to reverse direction as follows,
$$
\log(e^{is}-e^{i(t-i0)}) = i\pi + \log(e^{it} - e^{i(s+i0)})
 = i(\pi + t) + g(t-s)~,
$$
where $g$ is as in \eqref{g}.
Thus the exterior version of \eqref{vmspl} is
\be
v^+(Z(t)) := \frac{1}{2\pi} \int_{0}^{2\pi}
\left(\log\frac{e^{is}-e^{it}}{Z(s)-Z(t)}\right) \ttau(s) |Z'(s)| ds
+ \frac{T}{2\pi i}(\pi+t)
- \frac{1}{2\pi} \int_{0}^{2\pi}
g(t-s) \ttau(s) |Z'(s)| ds
\label{vpspl}
\ee
where the {\em total charge} of the single-layer density $\tau$ is
\be
T := \int_\pO\tau(y) |dy| = \int_{0}^{2\pi} \ttau(s) |Z'(s)| ds
~.
\label{T}
\ee
The first (smooth) term in \eqref{vpspl} is identical to that in \eqref{vmspl}.
The $\pi$ in the middle term can be dropped since it has no effect
on $\re v$.
In the last term, 
the convolution kernel is $g(-s)$ instead of $g(s)$;
we can achieve this by replacing $R_j$ with $R_{-j}$.
Thus the exterior version of \eqref{vmnyst} is
\be
v^+_k \approx
\frac{1}{2\pi} \sum_{j\neq k} \left( \log \frac{e^{i(s_k-s_j)}}{Z(s_k)-Z(s_j)}
\right) w_j \tau_j + \frac{1}{2\pi} \log \frac{i e^{is_k}}{Z'(s_k)} w_k \tau_k
+ \frac{T}{2\pi i} t_k
- \frac{1}{2\pi} \sum_{j=1}^N R_{k-j} w_j \tau_j
~.
\label{vpnyst}
\ee
This completes Step 1 for the Laplace single-layer exterior case.

\begin{rmk}
Kress \cite{kress91} gives formulae (due to Martensen--Kussmaul)
for splitting a periodic kernel with a logarithmic singularity
into a smooth part and the product of a smooth part and
the convolution kernel $\log \left (4 \sin^2 \frac{s}{2}\right)$.
The formula
$g(s) = \log(1-e^{-i(s-i0)}) = \half \log \left (4 \sin^2 \frac{s}{2}\right)
+ i\frac{s-\pi}{2}$, for $0\le s< 2\pi$,
shows that \eqref{split} is the analogous
Kress-type split for the {\em complex} logarithmic kernel case.
The imaginary part of $g$ is a ``sawtooth wave'' (periodized linear function)
whose sign depends on from which side the limit is taken.
\end{rmk}

If the total charge $T=0$
then $v$ given by \eqref{slpc} is single-valued in $\Oc$
and has $v_\infty=0$,
so Step 2 may proceed just as for the double-layer case, and we are done.
However, if $T\neq 0$, then there must be a branch cut in $v$
connecting $\pO$ to $\infty$
(along which the imaginary part must jump by $T$), and $v$
grows logarithmically at $\infty$.
In this latter case we must subtract off the total charge as follows.
Let $a\in\Omega$ be chosen not close to $\pO$, then define
\be
w(x)\; := \;v(x) - \frac{T}{2\pi}\log\frac{1}{a-x}~,
\label{w}
\ee
where $v$ is as in \eqref{slpc}.
Then $w$, being holomorphic in $\Oc$ (the branch cut of log in
\eqref{w} can be chosen to cancel that of $v$),
and having $w_\infty=0$,
is appropriate for representation by Step 2,
using its boundary data $w^+_k = v^+_k - (T/2\pi)\log 1/(a-y_k)$,
for $k=1,\ldots, N$.
Finally, after Step 2 produces $w$ and $w'$,
the missing monopole must be added back in
to get $v$ via \eqref{w}, and $v'(x) = w'(x) + (T/2\pi)/(a-x)$.
Then as before, $u = \re v$ and $\nabla u = (\re v', -\im v')$.

\begin{rmk}
The reader may wonder how spectral accuracy is to be achieved in Step 2
using the periodic trapezoid rule
when $T\neq 0$, given
that then the boundary data \eqref{vpnyst} has a {\em discontinuity}
due to the imaginary sawtooth (third term).
In fact there is also a discontinuity of equal size and opposite sign,
in $w$, introduced by the branch cut of the log in \eqref{w}.
Moreover, these discontinuities may occur at different boundary locations;
but, since the number of nodes lying between them is fixed,
the error introduced is a purely imaginary constant and has no effect on $u$.
The overall scheme for evaluating $u$ and $\nabla u$,
as for the interior case, is spectrally accurate.
\end{rmk}

\subsection{Numerical tests of Laplace evaluation quadratures}
\label{s:laptest}

In order to test the new global quadrature evaluation schemes
at a variety of distances from the boundary,
we set up simple BVPs on a smooth closed curve
with boundary data corresponding to known Laplace
solutions. We then solve each using an integral equation formulation
and record the error between the numerical layer potential evaluation and
the known solution.
Introducing and solving the BVP is necessary since on a general curve there
are very few known density functions which generate known analytic potentials
(specifically, $\tau\equiv1$ generating $u\equiv-1$ via the DLP
is the only example known to the authors).

The Laplace BVPs we use for tests are the standard four possibilities of
interior/exterior, Dirichlet/Neumann
problems \cite[Sec.~6.2]{LIE}. The integral equation
representations are those of \cite[Sec.~6.4]{LIE}
with the exception of the exterior Dirichlet, which we do not modify
since we choose the solution to vanish at $\infty$.
We use the representation $u={\cal D}\tau$
for the Dirichlet BVP, and $u={\cal S}\tau$ for the Neumann BVP,
with boundary data $f$.
Using the operator $D$ and the jump relations \eqref{sjr}--\eqref{djr},
the integral equations are thus,
\bea
\mbox{Dirichlet BVP}: \qquad & &(D \pm \half)\tau = f~,
\qquad (+ \mbox{ exterior case}, - \mbox{ interior case})
\label{lapdirbie}
\\
\mbox{Neumann BVP}: \qquad & & (D^T \mp \half)\tau = f~,
\qquad (- \mbox{ exterior case}, + \mbox{ interior case}) ~.
\label{lapneubie}
\eea
Since the kernel of $D$ is smooth, we
fill the system matrix via the Nystr\"om method \cite[Sec.~12.2]{LIE}
using the diagonal values $\lim_{s\to t}D(s,t) = -\kappa(t)/4\pi$
where $\kappa(t)$ is the curvature of $\pO$ at $Z(t)$.
The interior Neumann solution is only defined up to a constant,
hence we fix this constant by defining the error at the origin to be zero.
For a solution in $\Omega$ we use $u(x) = \re e^{i(1+x)}$,
$x\in\mathbb{C}$,
and in $\Oc$ we use $u(x) = \re 1/(x-0.1-0.3i)$,
which has $u_\infty=0$.
For the exterior Dirichlet and interior Neumann cases the operator,
and hence Nystr\"om matrix, has nullity 1; however, this does not pose
a problem when a backward-stable dense linear solver is used
(we use the backslash command in \matlab).

\bfi[t] 
\hspace{-2ex}\mbox{
\ig{width=1.9in}{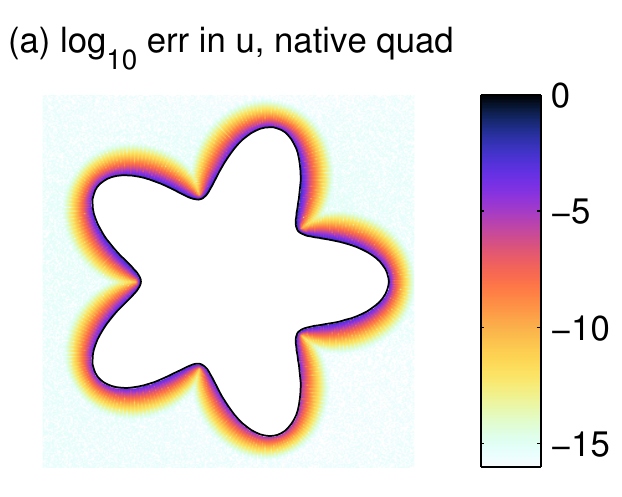}
\ig{width=2in}{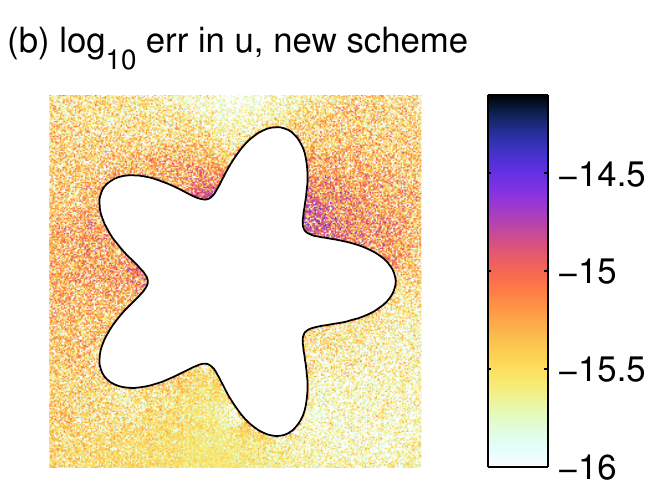}
\ig{width=2in}{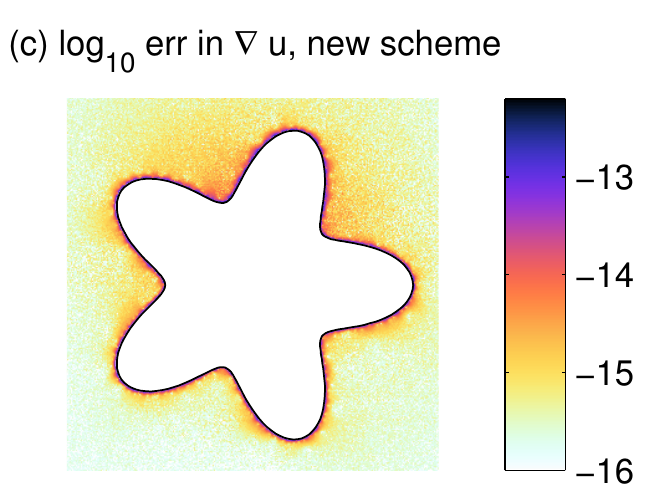}
}
\ca{
Evaluation error for the Laplace single-layer potential
on a grid exterior to the closed curve $\pO$ defined by the
radial function $r(\theta) = 1 + 0.3\cos 5\theta$.
The color shows $\log_{10}$ of the errors relative to a known
solution to the exterior Neumann BVP.
(a) uses the native periodic trapezoid rule \eqref{ptr}, whereas
(b) and (c) use the scheme of Section~\ref{s:lapslpext}
(note the change in color scale).
In all cases $N=240$ nodes are used for solution and evaluation.}{f:lapslpext}
\efi

\begin{table}[t] 
\centering
\begin{tabular}{l|rr|rr|rr|rr}
$N$ & \multicolumn{2}{c}{DLP int} & \multicolumn{2}{c}{DLP ext} & \multicolumn{2}{c}{SLP int} & \multicolumn{2}{c}{SLP ext}
\\
\hline
& $u$ & $\nabla u$ &$u$ & $\nabla u$ &$u$ & $\nabla u$ &$u$ & $\nabla u$
\\
\hline
100 & 2.9e-07&9.6e-06 & 8e-05 &2.6e-03  &   7e-09 &2.7e-07 & 1e-06 &3.9e-05\\
150 & 7.8e-11 &3.8e-09  & 6.7e-10 &6.8e-08    &  1.4e-12& 8.7e-11 & 7.9e-10 &7.5e-08 \\
200 & 2.1e-14 &2e-12 &  2.6e-13 &3.4e-11    &  9.8e-15 &7e-13     & 2.7e-13 &3.6e-11 \\
250 & 2e-14 &1.7e-12  & 4.7e-14 &4.6e-12    & 5.9e-14 &4.5e-12    & 4.9e-15 &6.3e-13
\\
\hline
\end{tabular}
\vspace{1ex}
\caption{Convergence of error for
the Laplace layer potential evaluation scheme
of Section~\ref{s:lap} for BVPs solved on the curve shown in
Figure~\ref{f:lapslpext}. The maximum error in $u$ or $\nabla u$
is taken over a square grid of spacing 0.01.
\label{t:lapconv}}
\end{table}

Figure~\ref{f:lapslpext} shows
results for the exterior single-layer potential
on the curve used in Sec.~\ref{s:cautest},
for the converged value of $N=240$.
As described in the introduction, we see that with
the native quadrature scheme
errors grow exponentially up to ${\mathcal O}(1)$ near $\Gamma$.
In contrast, the new scheme of this section achieves 14 digits in value,
and 12 digits in first derivative, uniformly for points
arbitrarily close to $\pO$.
The results for the other three BVPs are very similar.
For the convergence in all four cases we refer to Table~\ref{t:lapconv}
which shows the worst-case error over target points lying on a
square grid covering $[-1.5,1.5]^2$ with grid spacing 0.01
(some of these target points lie exactly on nodes; the closest
other ones are distance $3\times 10^{-4}$ from a node).
High-order convergence is apparent from this table, with the interior
schemes converging slightly faster in $N$ than the exterior ones.

\begin{rmk} 
Our test curve is chosen to be the same as the  interior Dirichlet BVP
test in \cite[Fig.~2--4]{helsing_close},
enabling a comparison of our scheme against their panel-based version.
When we use data from their solution,
$u(x) = \re 1/(z-1.5-1.5i) + 1/(z+0.25-1.5i) +1/(z+.5+1.5i)$,
we achieve uniform
14-digit accuracy by $N=320$; from \cite[Fig.~4]{helsing_close}
the panel-based version requires $N=480$.
So the periodic trapezoid rule is a factor 1.5 times more efficient in terms
of unknowns, which is close to the expected factor $\pi/2$ \cite{halequad}.
\label{r:piover2}
\end{rmk}

We use \matlab\ (R2012a) for this implementation (and others in this paper),
and achieve around $10^7$ source-target pairs per second on a 2.6 GHz i7 laptop.

\section{Stokes close evaluation scheme and numerical results}
\label{s:stokes}

We are at last in a position to describe how we evaluate Stokes potentials,
given samples $\{\sigma_j\}_{j=1}^N$
of the vector density $\sigma$ at $N$ trapezoid rule
nodes $\{y_j\}_{j=1}^N$ defining a closed curve $\Gamma$.
Throughout we generate all other geometric data (normals, curvature,
etc) at the nodes using spectral differentiation via the $N$-point FFT.

For the SLP we use \eqref{slp_stokes_via_lps},
where for each of the three Laplace SLPs we use the
samples of the density ($y\cdot\sigma$, $\sigma_1$ or $\sigma_2$)
at the same $N$ nodes, 
fed into the new Laplace evaluation scheme
of \sref{s:lapslp} (which itself relies on \sref{s:cau} for Step 2).
For the DLP we similarly use \eqref{dlp_stokes_via_lps},
with the last three terms using samples at the given $N$ nodes
fed into the (Ioakimidis/Helsing) scheme of \sref{s:lapdlp}.
However, for the first term in \eqref{dlp_stokes_via_lps},
we need the complex-valued densities $\tau_1$ and $\tau_2$
in \eqref{dlp_stokes_bowei_term}; notice that
the scheme of \sref{s:lapdlp} works perfectly well when fed a complex
$\tau$, producing \eqref{dlpcau} rather than \eqref{dlp}.
We have found that this first term converges
slower than the others,%
\footnote{We suspect this is because the two appearances
of the normal function in \eqref{dlp_stokes_bowei_term} cause
a more oscillatory integrand.}
so to preserve overall accuracy
we use FFT interpolation to upsample $\sigma$
by a factor $\beta>1$
when computing $\tau_1$ and $\tau_2$.
Thus this first term (and its two Laplace DLP evaluations) are
done with densities sampled on $\beta N$ nodes.
We have found that $\beta=2.2$ is sufficient to recover similar
accuracy to the other terms.%
\footnote{We note that panel-based schemes also need upsampling from
16 to 32 nodes per panel for full accuracy \cite{helsing_close}.}
This increases the effective cost of the Stokes DLP from 5 to about 7 Laplace
DLPs.

We now test the performance of our Stokes close evaluation
scheme in settings
relevant for vesicle simulations. We consider four examples:
the first two study the effect of proximity of the target
and complexity of the geometry on the error.
In the last two, similar
to the Laplace case, we set up BVPs with boundary data corresponding
to known analytic solutions and compare against numerical solutions
obtained through integral equation solves.
The fourth example is special in that we also use the evaluation scheme to
{\em apply}
the Nystr\"om matrix in the solve, given a boundary with close-to-touching
components.

\begin{figure}[t] 
\centering
\includegraphics[width=\textwidth]{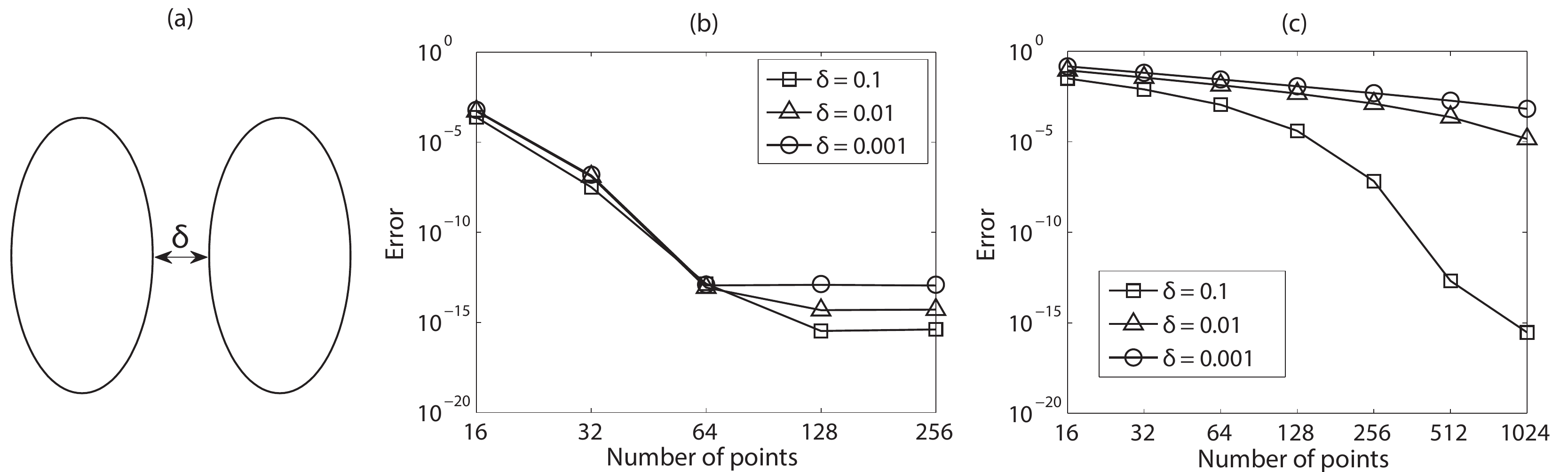}
\caption{Performance of proposed globally compensated Stokes SLP quadrature in
evaluating the hydrodynamic interaction force \eqref{interactionF}
between two elliptical bubbles parametrized as translations of
$(\cos \theta,2\sin \theta)$.
$\delta$ is the minimum distance between the two interfaces;
we plot the $l^{\infty}$-norm error against the number of nodes $N$ for $\delta = 0.1$, $0.01$ and $0.001$.
(b) Convergence of new globally compensated scheme.
(c) Native quadrature (trapezoidal rule), showing very poor convergence
for small $\delta$ values.
\label{error_ves_interact}}
\end{figure}

{\bf Example 1.} Consider two identical bubbles separated
by a small distance $\delta$ as shown in
Fig.~\ref{error_ves_interact}(a).  Assuming unit surface tension, the
interfacial force $f$ on each bubble is given by $f (y) = \kappa (y)
n_y$ where $\kappa$ is the curvature and $n_y$ is the unit outward
normal at a point $y$ on the interface $\Gamma$. The hydrodynamic interaction
force $F$ experienced at target $x$ on one bubble due to the presence
of the other is simply the Stokes SLP with interfacial force as the
density, that is,
\begin{equation} F(x) = (\SLP f) (x) = {1\over4\pi}\int\limits_{\Gamma}\left(\lir \, I + \dfrac{ r\otimes r}{\rho^2} \right) \kappa\, n_y \,ds_{y}
~. \label{interactionF}
\end{equation}
%
If $\delta$ is large, the integrand in \eqref{interactionF} is smooth
and the native trapezoidal rule will yield superalgebraic
convergence. When $\delta$ is small, this convergence rate is proportional
to $\delta$, due to the nearly singular integral;
hence the $N$ required scales like $1/\delta$, and is
unacceptably large as shown in
\fref{error_ves_interact}(c).
Yet our new scheme achieves 13 digits with only $N=64$ nodes, as shown in
\fref{error_ves_interact}(b).
%

\begin{figure}[t] 
\centering
\includegraphics[width=\textwidth]{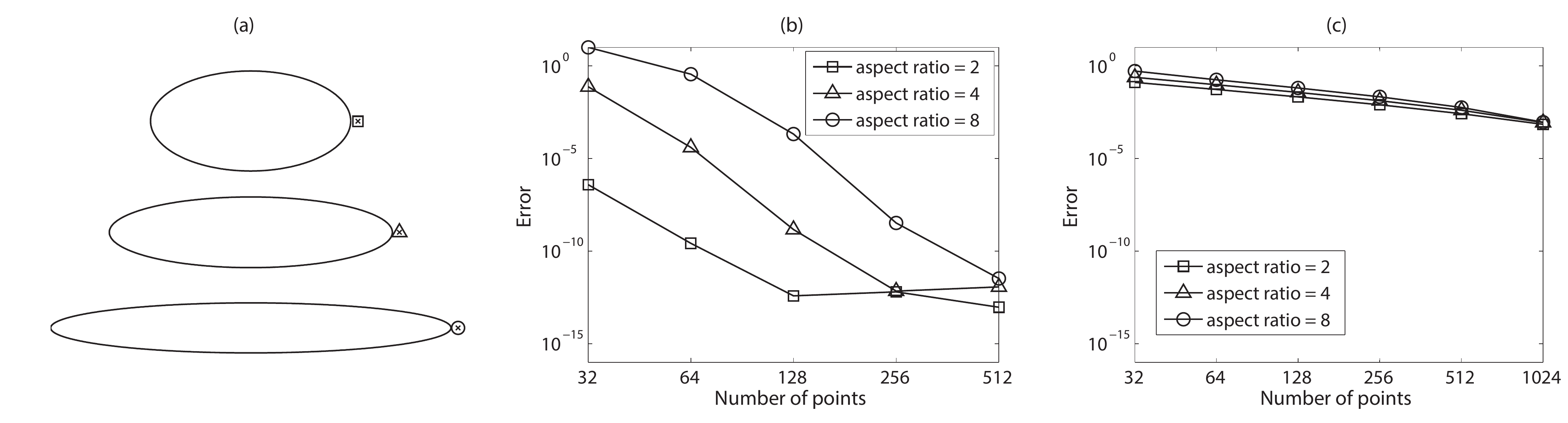}
\caption{Absolute errors in computing the Stokes SLP at a fixed distance $\delta=0.001$ away from the boundary of three different geometries with varying complexity: ellipses with aspect ratios $2, 4$, and $8$.
(a) Shows ellipses and target point $x$ (case $\delta=0.1$ is shown
to make the separation visible).
(b) Results for our proposed scheme.
In all cases, the errors decay exponentially with discretization size although, as expected, the absolute errors are higher for higher aspect ratio geometries. Therefore, in the case of our scheme, the given data (boundary, density) dictates the overall accuracy. (c) Results for native quadrature. In this case, convergence rate decreases but more striking is the fact the errors are almost the same for a particular number of points. This clearly implies that the near singular integral evaluation dominates the overall error.
\label{error_curvature}}
\end{figure}

{\bf Example 2.} We consider elliptical bubbles of aspect ratios 2, 4, and 8,
and evaluate $(\SLP f) (x)$ for a target point $x$ a distance $\delta=10^{-3}$
from the highest-curvature point of the ellipse.
This can be interpreted as the disturbance velocity at $x$
induced by the nearby bubble.
Fig.~\ref{error_curvature}(b) shows the rapid exponential
convergence of our SLP scheme.
For aspect ratio 2, by $N=128$ it has converged to 13 digits,
and for higher aspect ratios the $N$ required grows in proportion
to the aspect ratio.
This is as expected, since the density data ($\kappa n_y$) changes more
rapidly at the extreme point of the ellipse as its aspect ratio grows,
thus requires a larger $N$ to accurately interpolate.
Thus the errors are limited by resolving the data, not by the close-evaluation
scheme.
In contrast, errors using the native scheme are unacceptable,
and dominated by the small $\delta$.

\begin{figure}[t]  
\centering
\begin{tabular}{cc}
\includegraphics[height = 1.9 in]{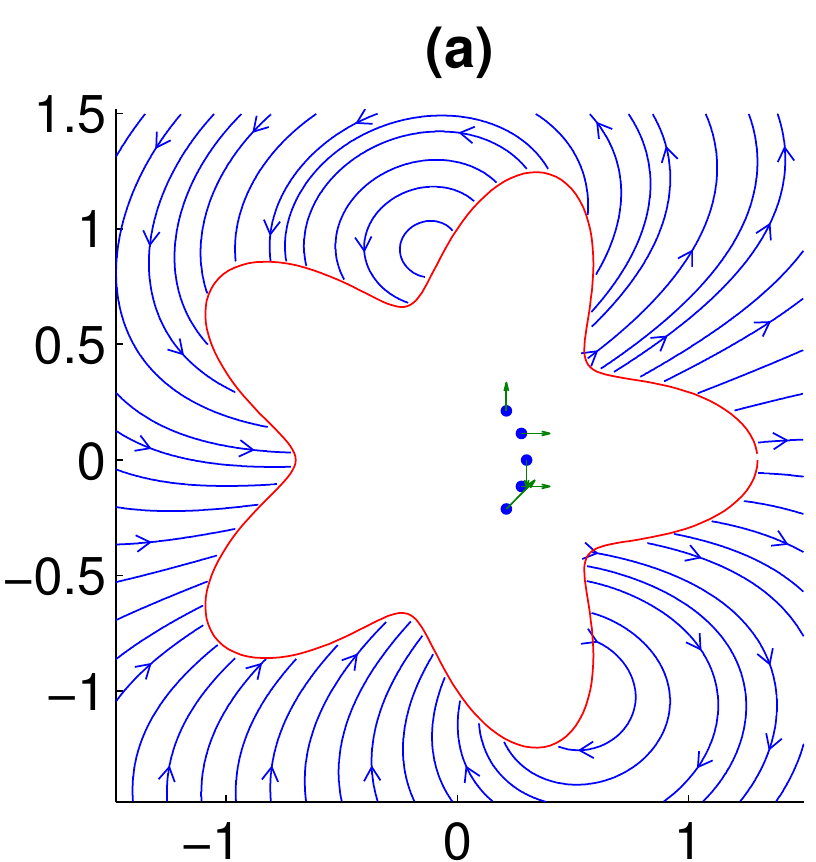}  &
\includegraphics[height = 1.9 in]{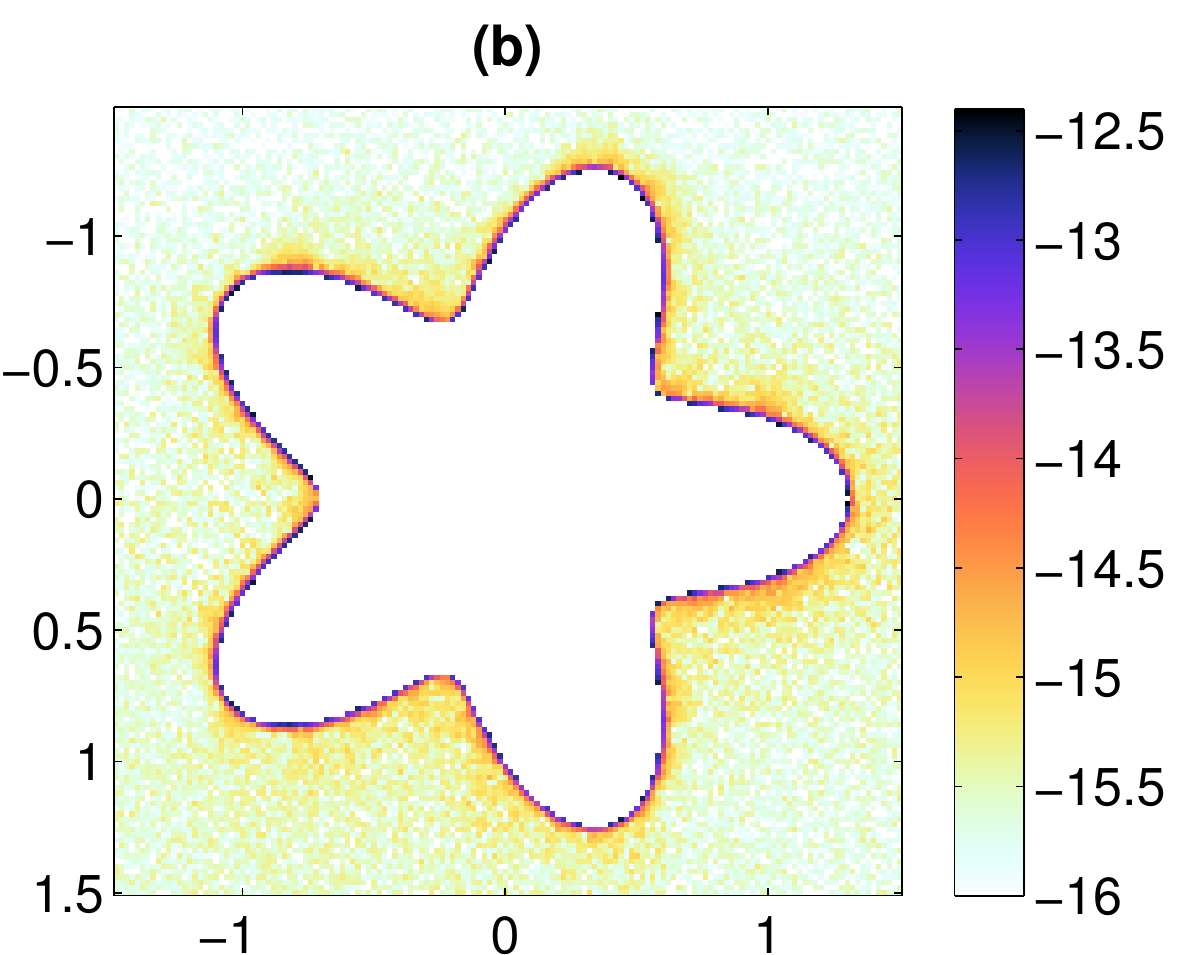}
\end{tabular}
\caption{(a) Streamlines of the velocity field generated by stokeslets shown in the interior, for testing the exterior BVPs in Example 3 of \sref{s:stokes}.
The boundary $\Gamma$ is given by the polar function $f(\theta) = 1 + 0.3 \cos 5\theta$.
(b) $\log_{10}$ of absolute error in computing the velocity
(on a grid of spacing 0.02) for the exterior Neumann Stokes
BVP, using $N = 250$ discretization points on the boundary for both the
Nystr\"om solve and the evaluation as in Example 3 of \sref{s:stokes}.
\label{error_ex3}}
\end{figure}

\begin{table}[ht]  
\centering
\begin{tabular}{c|l|l|l|l}
$N$ & ext Dirichlet & int Dirichlet & ext Neumann & int Neumann\\
\hline
100 & 2.3e-05 & 3.6e-05 & 2.7e-04 & 1.8e-04\\
150 & 2.3e-07 & 9.9e-08 & 6.1e-08 & 1.4e-08\\
200 & 5.0e-09 & 1.9e-09 & 3.7e-12 & 4.0e-11\\
250 & 2.8e-10 & 1.4e-11 & 3.6e-13 & 1.2e-12\\
300 & 2.3e-11 & 8.3e-14 & 6.7e-13 & 2.0e-13\\
350 & 2.0e-12 & 4.9e-14 & 4.3e-13 & 1.1e-13\\
\hline
\end{tabular}
\caption{Convergence of velocity error
(max over the evaluation grid) for the solutions of the four
types of Stokes BVPs, on the curve shown in \fref{error_ex3}.
The Nystr\"om method is used to find the density,
then the new close evaluation scheme is
used, as described in \sref{s:stokes}.}
\label{t:ex3}
\end{table}





{\bf Example 3.} We solve all four types (interior/exterior,
Dirichlet/Neumann) of Stokes BVP on
the star-shaped geometry introduced in Section \ref{s:cautest}
using a standard Nystr\"om method,
and test the velocity field evaluation error with the new scheme.
Denoting the pressure by $p$ and the boundary vector data by $g$,
the (unit viscosity) Stokes BVPs are \cite[Sec.~2.3.2]{HW}
\bea
-\Delta u + \nabla p &=& 0 \quad
\mbox{ in } \Oc
\mbox{ (exterior)~, \quad or in } \Omega \mbox{ (interior)~,} \\
\nabla \cdot u &=& 0 \quad
\mbox{ in } \Oc
\mbox{ (exterior)~, \quad or in } \Omega \mbox{ (interior)~,} \\
u &=& g \mbox{ on $\Gamma$ \; (Dirichlet), \quad or }
T(u,p) \; = \; g \mbox{ on $\Gamma$ \; (Neumann)}~,
\label{bc} \\
u(x) &=& \bm\Sigma \log \rho + {\mathcal O}(1) \qquad \mbox{ as }
\rho:=|x| \to\infty \quad \mbox{ (exterior only)~,}
\eea
where the traction at $y\in\Gamma$ is defined by $T(u,p):= -p n_y +
\bigl(\nabla u + (\nabla u)^T\bigr)\otimes n_y$, and
$\bm\Sigma$ is some constant vector.
For the Dirichlet case $g$ is the boundary velocity, for the Neumann
case the boundary traction.

We use standard representations as follows.
Let $\Sop$ and $\Dop$ denote the boundary integral operators with
the Stokes SLP and DLP kernels,
with $\Dop$ taken in the principal value sense
(these are analogous to $S$ and $D$ in the Laplace case).
Then for the Neumann BVPs we use $u = \SLP \sigma$,
giving the integral equation $(\Dop^T \mp \half)\sigma = g$
($-$ exterior case, $+$ interior case), just as in \eqref{lapneubie}.
For the interior Dirichlet BVP we use $u= \DLP \sigma$,
giving $(\Dop - \half)\sigma = g$.
For the exterior Dirichlet, since we wish to solve data $g$ for which
$\Sigma \neq \mbf{0}$, we use the combined representation
$u = (\DLP + \SLP) \sigma$
\cite[page 128]{pozrikidis},
giving $(\Dop +\Sop + \half)\sigma = g$.
%
%
This also eliminates the 1-dimensional nullspace in this case.%
\footnote{Traditionally, a combination of stokeslets and rotlets are used to eliminate this nullspace in the context of Stokes BVPs. A more general procedure for handling nullspaces is recently given in \cite{sifuentes2014randomized}.}
All four integral equations are consistent,
but the interior Dirichlet and exterior Neumann have
1-dimensional nullspaces which only affect the pressure solution $p$,
and the interior Neumann integral equation has a 3-dimensional
nullspace corresponding to rigid motions in the plane 
\cite[Sec.~2.3.2]{HW}.
In the last case we match these three components to those of the
reference solution before computing the errors.

We use the boundary $\Gamma$ shown in \fref{error_ex3} for all four BVP
tests.
A simple exterior reference solution $(u^*,p^*)$ is constructed via
stokeslets placed at random locations in the interior. At any exterior point $x$, the velocity and pressure are then given by
\bea
u^*(x) &\;=\;&
\frac{1}{4\pi}\sum_{k=1}^q \biggl(\log \frac{1}{|x - y_k|} \biggr)f_k +
\dfrac{ (x - y_k)\cdot f_k}{|x-y_k|^2} \, (x - y_k) \\
p^*(x) &\;=\;&
\frac{1}{2\pi}\sum_{k=1}^q \frac{(x-y_k)\cdot f_k}{|x-y_k|^2}
\label{sol}
\eea
The number of stokeslets $q$ is set to 5 and their
chosen locations $\{y_k\}$ and strengths $\{f_k\}$
are displayed in Fig.~\ref{error_ex3}(a).
For the interior cases the stokeslets are moved to the exterior of $\Gamma$
at a radius of 2. 
Then either the Dirichlet data $g(x) = u^*(x), x\in\pO$
or the Neumann data $g(x) = T(u^*,p^*)(x), x\in\pO$ is used.
The sup norm of this data is always in the range $0.1$ to 1.

The numerical velocity solution $u(x)$ is obtained by solving the
integral equation for $\sigma$ using the standard Nystr\"om method,
then using the new Stokes close evaluation schemes for the SLP and/or DLP
to evaluate the representation for $u$.
We note that for the exterior Dirichlet case, the logarithmically-singular
operator $\Sop$ must be discretized; for this we use 16th-order Alpert
corrections \cite{alpert} as explained in \cite[Sec.~4]{hao}.
For the exterior Neumann case, the resulting error is plotted on a grid
of spacing 0.02 in \fref{error_ex3}(b).
The convergence of worst-case errors on this grid for all four BVPs is
given in \tref{t:ex3}.
It is superalgebraic, and clear that
12-digit accuracy results for $N$ in the range 250 to 300,
apart from the exterior Dirichlet case when 350 is needed.%
\footnote{We verified that it was the Alpert correction for $\Sop$ that caused
this slightly slower convergence rate.}
Away from a thin boundary layer the errors improve by at least 2 digits.
Note that it is not very meaningful to do a direct comparison against
the Laplace case since the smoothness of the data
(i.e.\ distance away of the sources $y_k$) has a large effect
on the convergence rate.

\bfi  
\centering\ig{width=\textwidth}{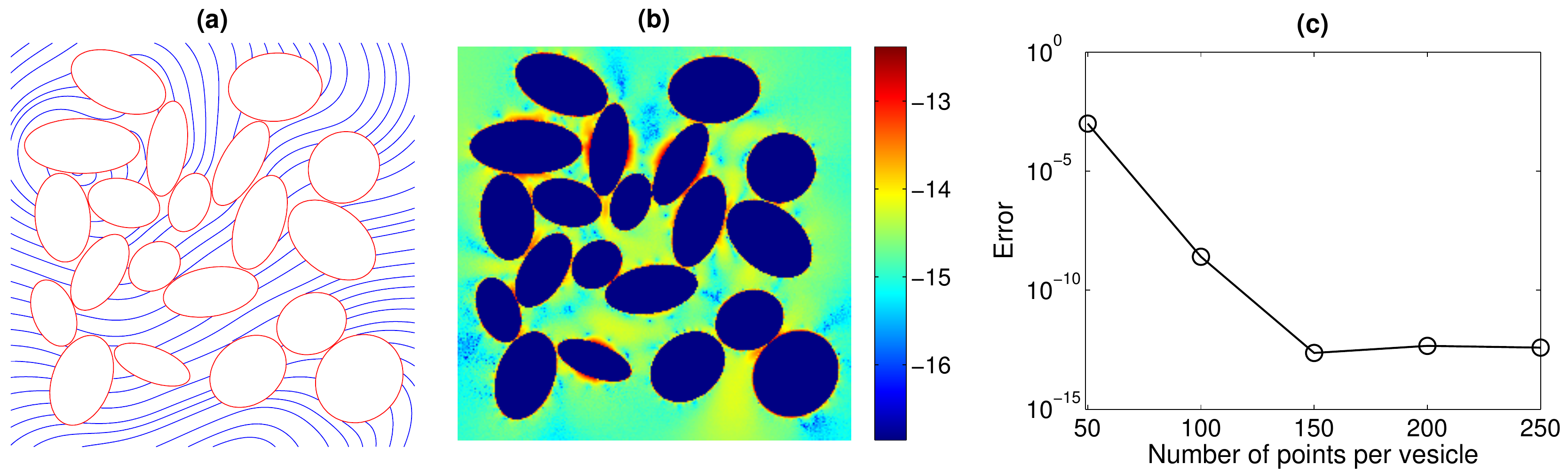}
\ca{Test of exterior Dirichlet BVP for the 20 elliptical vesicles of Example 4,
using the scheme of \sref{s:stokes}.
(a) Streamlines of the velocity field in the exterior of vesicles.
(b) $\log_{10}$ of error in velocity field using
the Nystr\"om method (with GMRES) to solve the
density, using the close evaluation scheme for vesicle-vesicle
interactions and for final flow evaluation, with $N=150$ nodes per vesicle.
(c) convergence of the sup norm of error in the velocity field
over the grid of values plotted in (b).
}{20ves}
\efi

{\bf Example 4.}
Our final test focuses on the exterior Dirichlet Stokes BVP, using
the same combined representation $u=(\DLP +\SLP)\sigma$ as before,
with a more complicated boundary $\Gamma$ comprising 20 elliptical
vesicles that come very close to each other, shown in \fref{20ves}(a).
Their minimum separation is $2\times 10^{-4}$.
This geometry is taken from \cite[Table~4]{yingbeale}.
The number of quadrature nodes for each vesicle is the same, $N$.
To resolve the Nystr\"om matrix for the solution of $\sigma$
one would need very high $N$ when using the standard Nystr\"om
quadrature formulae, due to interacting vesicles.
Therefore in this example we use an iterative solver (GMRES),
using the new global close evaluation scheme for the SLP and DLP
to apply all vesicle-vesicle interaction blocks of the system matrix.
(The self-interactions of each vesicle are done by the standard
Nystr\"om formulae, with 8th-order Alpert corrections.)
The number of GMRES iterations needed was high (around 600);
however,
our point is merely to demonstrate that our close evaluation scheme can be
used to solve a close-to-touching geometry with a small $N$ per vesicle.
%

The reference solution is the flow field due to a single stokeslet
of random strength at the center of each ellipse.
We believe that this test case is crudely representative of the types
of vesicle flows occurring in applications, with
non-singular velocities where vesicles approach each other (this contrasts
the case of approaching rigid bodies).
The solution error for $u$ is evaluated on a grid of spacing 0.016
for $N=150$ in \fref{20ves}(b), and the convergence
tested with respect to $N$ in  \fref{20ves}(c).
Convergence is again superalgebraic, reaching 13 digits by $N=150$.
By comparison, the recent tests of \cite{yingbeale} achieve only around 3 digits
at this $N$.

\section{Conclusions and discussion}
\label{s:conclusion}

We have presented a simple new scheme
for evaluating the classical Stokes layer potentials on smooth closed curves
with rapid spectral convergence in the number of quadrature
nodes $N$, independent of the distance from the curve.
We expect this to find applications in Stokes flow simulations
with large numbers of close-to-touching vesicles.
This builds upon Laplace
layer potentials---including a new scheme for
single-layer potentials of independent interest---which
%
%
in turn rely on barycentric-type quadratures for Cauchy's formula.
The underlying global periodic trapezoid rule,
being simpler and potentially more efficient than panel-based schemes
(see \rref{r:piover2}), is the most common in vesicle applications.
Thus our work complements recent panel-based Stokes quadratures
\cite{ojalastokes}.


We performed systematic tests, solving all eight boundary value
problems (Laplace/Stokes, Dirichlet/Neumann, interior/exterior)
via integral equations,
and reach close to machine error in the solution with a small $N$ that is
essentially the same as that needed for the Nystr\"om method itself.
%
We also showed that our evaluator can be used
to apply the operator in an iterative solution with close-to-touching
boundaries (Example 4).

To evaluate at $M$ targets the effort is ${\mathcal O}(N(N+M))$,
although this would be easy to improve to ${\mathcal O}(N+M)$,
%
%
by using the Cauchy FMM in the barycentric
evaluations \eqref{vcomp} and \eqref{vcompe} and the Laplace Step 1
(splitting the sums in \eqref{vmk}--\eqref{vpk}),
and a complex logarithmic FMM in \eqref{vmnyst} and \eqref{vpnyst}.
In vesicle applications $N$ is small, so we leave this for future work.



In terms of future work, there are several variants that could be benchmarked,
such as whether derivatives are best moved from Step 2 to Step 1,
and alternatives in Remarks~\ref{r:2pii} and \ref{r:tref}.
However, given the results we presented, there is not
much room for improvement in accuracy.
For optimal speed, a Fortran/C/OpenMP library is certainly needed.
%
Since the Cauchy formula quadratures of \sref{s:cau}
(due to Ioakimidis) generalize barycentric interpolation to the complex
plane, this connection is worth analyzing further
(along these lines see \cite{austinrou}).
Finally, we do not know of any way to extend the
complex-analytic methods presented to 3D. For surfaces in 3D
one cannot yet say which of the many existing schemes is to be preferred.

We provide documented MATLAB codes
for the close evaluation of 2D Laplace and Stokes layer potentials,
and driver codes to generate the tables and several figures from this paper,
at the following URL:

{\tt http://math.dartmouth.edu/$\sim$ahb/software/lsc2d.tgz}

\section*{Acknowledgments}
We benefited from discussion with Nick Trefethen and Gary Marple.
This work of AHB was supported by the National Science Foundation
under grant DMS-1216656.  The work of SKV was supported by the
National Science Foundation under grant DMS-1224656.



\bibliography{alex}
\end{document}